\documentclass[a4paper,10pt]{article}
\title{Representation theory of super Yang-Mills algebras}
\author{Estanislao Herscovich 
\thanks{The author is an Alexander von Humboldt research fellow.}}
\date{}

\usepackage{amsmath,amsthm,amsfonts,amssymb,stmaryrd}
\usepackage{latexsym}
\usepackage{color}
\usepackage{graphicx}
\usepackage{float}  
\usepackage{fullpage}
\usepackage[small, bf, margin=90pt, tableposition=bottom]{caption}
\usepackage[T1]{fontenc}
\usepackage{palatino}
\usepackage[fleqn,tbtags]{mathtools}

\usepackage[
  breaklinks,
  naturalnames,
  ]{hyperref}

\usepackage{amsrefs}

\input{xy}
\xyoption{all}
\xyoption{poly}
\usepackage[all]{xy}

\newtheorem{theorem}{Theorem}[section]
\newtheorem{theorem*}{Theorem*}
\newtheorem{corollary}[theorem]{Corollary}

\newtheorem{proposition}[theorem]{Proposition}
\newtheorem{remark}[theorem]{Remark}

\numberwithin{equation}{section}                    
\def\id{{\mathrm{id}}}
\let\oldqed\qed
\renewcommand\qed{\oldqed\par\bigskip}

\newcommand\cl[1]{{\langle#1\rangle}}

\newcommand\CC{{\mathbb{C}}}
\newcommand\RR{{\mathbb{R}}}

\newcommand\ZZ{{\mathbb{Z}}}
\newcommand\NN{{\mathbb{N}}}

\newcommand\gl{{\mathfrak{gl}}}

\newcommand\so{{\mathfrak{so}}}

\newcommand\SO{{\mathrm{SO}}}
\def\Aut{{\mathrm {Aut}}}
\def\Prim{{\mathrm{Prim}}}

\def\B{{\mathcal B}}
\def\C{{\mathcal C}}

\def\U{{\mathcal U}}

\def\YM{{\mathrm {YM}}}
\def\TYM{{\mathrm {TYM}}}
\def\Cliff{{\mathrm {Cliff}}}

\def\ad{{\mathrm {ad}}}

\def\Ker{{\mathrm {Ker}}}
\def\Hom{{\mathrm {Hom}}}
\def\Mod{{\mathrm {Mod}}}
\def\mod{{\mathrm {mod}}}

\def\End{{\mathrm {End}}}

\def\Der{{\mathrm {Der}}}
\def\InnDer{{\mathrm {InnDer}}}

\def\Tor{{\mathrm {Tor}}}

\def\a{{\mathfrak a}}
\def\g{{\mathfrak g}}
\def\h{{\mathfrak h}}
\def\k{{\mathfrak k}}

\def\ym{{\mathfrak{ym}}}
\def\f{{\mathfrak f}}

\def\tym{{\mathfrak{tym}}}

\def\place{{-}}

\def\MyNode{\ifcase\xypolynode\or
      (W \otimes X) \otimes (Y \otimes Z)
    \or
      ((W \otimes X) \otimes Y) \otimes Z
    \or
      (W \otimes (X \otimes Y)) \otimes Z
    \or
      W \otimes ((X \otimes Y) \otimes Z)
    \or
      W \otimes (X \otimes (Y \otimes Z))
    \fi
  }%

\def\MyNodes{\ifcase\xypolynode
    \or
      X \otimes Y
    \or
      X \otimes (e \otimes Y)
    \or
     (X \otimes e) \otimes Y
    \fi
  }%

\def\MyNodessr{\ifcase\xypolynode
    \or
      X \otimes Y
    \or
      e \otimes (X \otimes Y)
    \or
     (e \otimes X) \otimes Y
    \fi
  }%

\def\MyNodessl{\ifcase\xypolynode
    \or
      X \otimes Y
    \or
      (X \otimes Y) \otimes e
    \or
       X \otimes (Y \otimes e)
    \fi
  }%
\begin{document}

\maketitle
\begin{abstract}
   {
    We study in this article the representation theory of a family of super algebras, called the \emph{super Yang-Mills algebras}, 
    by exploiting the Kirillov orbit method \textit{\`a la Dixmier} for nilpotent super Lie algebras. 
    These super algebras are a generalization of the so-called \emph{Yang-Mills algebras}, introduced by A. Connes and M. Dubois-Violette in \cite{CD02}, 
    but in fact they appear as a ``background independent'' formulation of supersymmetric gauge theory considered in physics, 
    in a similar way as Yang-Mills algebras do the same for the usual gauge theory. 
    Our main result states that, under certain hypotheses, all Clifford-Weyl super algebras $\Cliff_{q}(k) \otimes A_{p}(k)$, 
    for $p \geq 3$, or $p = 2$ and $q \geq 2$, appear as a quotient of all super Yang-Mills algebras, for $n \geq 3$ and $s \geq 1$. 
    This provides thus a family of representations of the super Yang-Mills algebras. 
   }
\end{abstract}



\section{Introduction}

This article is devoted to the study of the representation theory of super Yang-Mills algebras.
Let us briefly recall the definition of the super Yang-Mills algebras. 
Given two nonnegative integers $n, s \in \NN_{0}^{2} \setminus \{ (0,0) \}$, and a collection of $(s \times s)$-matrices $\Gamma_{a,b}^{i}$, 
for $i = 1, \dots, n$ ($a,b = 1, \dots, s$), 
the \emph{super Yang-Mills algebra} $\ym(n,s)^{\Gamma}$ over an algebraically closed field $k$ of characteristic zero is defined as 
the quotient of the free super Lie algebra $\f(x_{1},\dots,x_{n},z_{1},\dots,z_{s})$, for even indeterminates $x_{1}, \dots, x_{n}$ and odd ones 
$z_{1}, \dots, z_{s}$, by the (homogeneous) relations given by  
\begin{align*}
 r_{0,i} &= \sum_{j=1}^{n} [x_{j},[x_{j},x_{i}]] - \frac{1}{2} \sum_{a,b = 1}^{s} \Gamma^{i}_{a,b} [z_{a},z_{b}], 
 \\
 r_{1,a} &= \sum_{i=1}^{n} \sum_{b = 1}^{s} \Gamma_{a,b}^{i} [x_{i},z_{b}],  
\end{align*} 
for $i = 1, \dots, n$ and $a = 1, \dots, s$, respectively. 
It can also be regarded as a graded Lie algebra with $\deg(x_{i})=2$, for $i = 1, \dots, n$, and $\deg(z_{a})=3$, for $a = 1, \dots, s$. 
The associative enveloping algebra $\U(\ym(n,s)^{\Gamma})$ will be denoted $\YM(n,s)^{\Gamma}$. 
They have been previously considered by M. Movshev and A. Schwarz in \cite{Mov05} and \cite{MS06}. 
Also, the case $s = 0$ (and $n \geq 2$) leads to the definition of Yang-Mills algebra given by A. Connes and M. Dubois-Violette in \cite{CD02}. 
Omitting the trivial cases with $n=0$, or the already known ones with $s=0$, 
we shall see that, for each $(n,s) \in \NN^{2}$, it is noetherian if and only if $n = 1$ and $s = 1, 2$. 
However, it is coherent for all values of the parameters $(n,s)$ (see Remark \ref{rem:coh}). 

From the physical point of view, it can be seen that the components of the covariant derivative and the dual spinor field of a supersymmetric gauge theory 
on the Minkowski space provide a representation of the corresponding super Yang-Mills algebra (\textit{cf.} \cite{DF99s}, and see also Remark \ref{rem:phys}). 
Otherwise stated, these super algebras yield a ``background independent'' formulation of the supersymmetric gauge theory in physics. 
They can also be regarded in order to provide a noncommutative version of it. 
Our interest in them comes then in order to shed more light in this direction. 

The main result of this article may be formulated as follows:
\begin{theorem*}
Let $n, s, p, q \in \NN$ be positive integers, satisfying $n \geq 3$. 
We suppose further that either $p \geq 3$, or $p = 2$ and $q \geq 2$. 
Then, there exists a surjective homomorphism of super algebras
\[       \U(\ym(n,s)^{\Gamma}) \twoheadrightarrow \Cliff_{q}(k) \otimes A_{p}(k),     \]
where $\Cliff_{q}(k) \otimes A_{p}(k)$ denotes the Clifford-Weyl super algebra. 
Furthermore, there exists $l \in \NN$ such that we can choose this morphism in such a way that it factors through the quotient
$\U(\ym(n,s)^{\Gamma}/F^{l}(\ym(n,s)^{\Gamma}))$, where $F^{l}(\ym(n,s)^{\Gamma})$ is the Lie ideal of $\ym(n,s)^{\Gamma}$ formed by the elements 
of degree greater than or equal to $l+2$.
\end{theorem*}

We would like to remark that the Clifford-Weyl super algebra $\Cliff_{q}(k) \otimes A_{p}(k)$ appearing in the previous theorem has the $\ZZ/2\ZZ$-grading given by usual grading of the Clifford (super) algebra 
$\Cliff_{q}(k)$ and by considering the Weyl algebra $A_{p}(k)$ to be concentrated in degree zero (see \cite{Her10}, Example 1.2). 
This differs from the grading of the ``Clifford-Weyl algebras'' $\mathcal{C}(q, 2 p)$ considered in \cite{MPU09}, since in that case the Weyl algebra has also a nontrivial homogeneous component of odd degree. 

In order to prove the theorem we needed to extend the so called \emph{Dixmier map} for nilpotent Lie algebras to the case of nilpotent super Lie algebras, 
which was done in \cite{Her10}. 
Even tough this article may be considered as an extension and generalization of the results proved in \cite{HS10}, 
it also deals with several difficulties and differences with the latter, which mainly follow from the fact that the study of the 
enveloping algebras of super Lie algebra has various differences with the case of enveloping algebras of Lie algebras. 
To mention just a few, the super Yang-Mills algebras are not Koszul, at least not for any definition we are aware of, 
even though they behave quite the same, 
there are not any a priori morphism between different super Yang-Mills algebras for arbitrary $\Gamma$'s (see the paragraph at the end of Subsection \ref{subsec:defi}), 
there are several differences between the Dixmier map of nilpotent Lie algebras and of nilpotent super Lie algebras 
(\textit{e.g.} the super dimension of a polarization at an even functional does not determine the weight of the ideal that it defines, which explains the phenomena that appears in Remark \ref{rem:otrosn}, etc), 
the enveloping algebra of a super Lie algebra is not necessarily semiprimitive, and the determination of its radical is usually a highly nontrivial task, etc. 

\bigskip

The contents of the article are as follows.
In Section \ref{sec:generalities} we recall the definition and several elementary properties of super Yang-Mills algebras, 
some of them with a physical flavour. 
In particular, we show in Subsection \ref{subsec:pot} that the typical odd supersymmetries of the supersymmetric gauge theories considered in physics also appear, under the same assumptions as there, 
as super derivations of the super Yang-Mills algebra. 

In Section \ref{sec:hom} we study the homological properties of this family of super algebras. 
In fact, Subsection \ref{subsec:projres} provides a complete description of the minimal projective resolution of the trivial module $k$ 
over the graded Lie algebra $\ym(n,s)^{\Gamma}$, for $(n,s) \neq (1,0),(1,1)$, 
which satisfy a property similar to koszulity. 
Using a procedure similar to the one employed by R. Berger and N. Marconnet in \cite{BM06}, Section 4, we furthermore obtain 
the minimal projective resolution of the super Yang-Mills algebra $\YM(n,s)^{\Gamma}$, considered as a graded algebra, in the category of bimodules, 
for the same set of indices. 
From the particular description of these minimal projective resolutions we prove that $\YM(n,s)^{\Gamma}$ AS-regular in the sense of \cite{MM11} 
and graded Calabi-Yau, for $(n,s) \neq (1,0), (1,1)$.   
We later derive some consequences, computing in particular the Hilbert series of both $\ym(n,s)^{\Gamma}$ and $\YM(n,s)^{\Gamma}$. 
Moreover, we prove that, for $n \geq 2$, or $n = 1$ and $s \geq 3$, the super Yang-Mills algebra $\ym(n,s)^{\Gamma}$ contains a finite codimensional Lie ideal 
which is a free super Lie algebra. 
This is done using simpler methods than the ones used in \cite{HS10}, Section 3. 
From this fact, we derive several consequences, and in particular that these graded algebras are not noetherian, but they are coherent. 

Finally, in section \ref{sec:main} we prove our main result, Theorem \ref{teo:yangmillsweyl},
and describe the families of representations appearing in this way.

We would like to thank A. Solotar for several suggestions and remarks. 

\section{Generalities}
\label{sec:generalities}

In this first section we fix notations and recall some elementary properties of what we call the super Yang-Mills algebras.

\subsection{Definition}
\label{subsec:defi}

Throughout this article $k$ will denote an algebraically closed field of characteristic zero.
The main convention and notations on super vector spaces, super algebras and modules over them (and also for the graded analogous ones) 
that we follow are the same as in \cite{Her10}, to which we refer. 
Unless otherwise stated, a module over an algebra (resp. a graded algebra, a super algebra), will always denote a left module. 
Moreover, we consider the category of modules over graded (resp. super) algebras provided with homogeneous linear morphisms of degree zero. 
It is also endowed with the shift functor $(\place)[1]$, defined by $(M[1])_{n} = M_{n+1}$, for $n \in \ZZ$ in the graded case, and 
$n \in \ZZ/2\ZZ$ in the super case. 
Given two modules $M$ and $N$ over a graded (resp. super) algebra $A$, $\mathrm{hom}_{A}(M,N)$ will stand for space of morphisms 
in the previously described categories. 
Furthermore, the internal space of morphisms is given by $\mathcal{H}om_{A}(M,N) = \oplus_{i \in G} \mathrm{hom}_{A}(M,N[i])$, where $G = \ZZ$ or $G = \ZZ/2\ZZ$ for the graded or the super case, respectively. 

We fix the following setup. 
Let $V = V_{0} \oplus V_{1}$ be a super vector space over $k$ of super dimension $(n,s) \in \NN_{0}^{2}$ with $n + s > 0$, 
such that the even part $V_{0}$ is provided with a nondegenerate symmetric bilinear form $g$. 
Note that the algebraic group $\SO(V_{0},g)$, and hence the Lie algebra $\so(V_{0},g)$, acts on $V_{0}$ (with the standard action). 
We shall also write $V = V(n,s)$, $V_{0}=V(n)_{0}$ or $V_{1}=V(s)_{1}$ if we want to stress the (super) dimension. 
We suppose further that there exists a map of the form $\Gamma : S^{2}V^{*}_{1} \rightarrow V_{0}$. 

Choose a (homogeneous) basis $\B = \B_{0} \cup \B_{1}$ of $V$, where $\B_{0} = \{ x_{1}, \dots, x_{n} \}$ and $\B_{1} =\{ z_{1}, \dots, z_{s} \}$, 
with $|x_{i}| = 0$, for all $i = 1, \dots, n$, and $|z_{a}| = 1$, for all $a = 1, \dots, s$, 
and let $\B^{*} = \B^{*}_{0} \cup \B^{*}_{1}$, where $\B^{*}_{0} = \{ x_{1}^{*}, \dots, x_{n}^{*}\}$ and $\B^{*}_{1} = \{z_{1}^{*}, \dots, z_{s}^{*} \}$, 
be the dual basis of $V^{*}$. 
Set $\Gamma_{a,b}^{i} = x_{i}^{*}(\Gamma(z_{a}^{*},z_{b}^{*}))$, for $i = 1, \dots, n$ and $a, b = 1, \dots, s$, and 
$g^{-1}$ the \emph{inverse} nondegenerate symmetric bilinear form on $V^{*}_{0}$, 
\textit{i.e.} $g^{-1}$ is the bilinear form on $V^{*}_{0}$ defined as the image of $g$ under 
the $k$-linear isomorphism $V_{0} \rightarrow V^{*}_{0}$ given by $v \mapsto g(v,\place)$. 
It is easy to see that the matrix of $g^{-1}$ with respect to the dual basis $\B^{*}_{0}$ is the inverse of the matrix of $g$ with respect to the basis $\B_{0}$, 
and we further write $g^{i,j} = g^{-1}(x_{i}^{*},x_{j}^{*})$ and $g_{i,j} = g(x_{i},x_{j})$. 

If $\f(V)$ denotes the free super Lie algebra generated by the super vector space $V$,
the \emph{super Yang-Mills algebra} is defined as the quotient 
\[     \ym(V,g)^{\Gamma} = \f(V)/\cl{R(V,g)^{\Gamma}},     \]
where $R(V,g)^{\Gamma}$ is the super vector space inside $\f(V)$ spanned by the elements 
\begin{equation}
\label{eq:rel}
\begin{split}
 r_{0,i} &= \sum_{j,l,m=1}^{n} g^{j,l} g^{i,m} [x_{j},[x_{l},x_{m}]] - \frac{1}{2} \sum_{a,b = 1}^{s} \Gamma^{i}_{a,b} [z_{a},z_{b}], 
 \\
 r_{1,a} &= \sum_{i=1}^{n} \sum_{b = 1}^{s} \Gamma_{a,b}^{i} [x_{i},z_{b}],  
\end{split}
\end{equation}
for $i = 1, \dots, n$ and $a = 1, \dots, s$, respectively. 

We consider the universal enveloping algebra $\YM(V,g)^{\Gamma} = \U(\ym(V,g)^{\Gamma})$ of $\ym(V,g)^{\Gamma}$, and also call it the 
\emph{(associative) super Yang-Mills algebra}. 
By definition, it is the super algebra given by the quotient of the tensor algebra $TV(n,s)$ of the super vector space $V(n,s)$ by the two-sided ideal generated by 
the same super vector subspace $R(V(n,s),g)^{\Gamma}$, now seen inside of $F^{3}TV(n,s)$, 
where $\{ F^{\bullet}TV(n,s) \}_{\bullet \in \NN_{0}}$ denotes the canonical (increasing) filtration of the tensor algebra. 
Note that $\U(\ym(n,s)^{\Gamma})$, being the enveloping algebra of a super Lie algebra, need not be a domain, whereas 
$\U(\ym(n,s)^{\Gamma}_{0})$ is always so. 
However, we shall see below that $\U(\ym(n,s)^{\Gamma})$ is also a domain. 

It is direct to check that $R(V,g)^{\Gamma}$ is independent of the choice of the homogeneous basis $\B$, so we may suppose that $\B_{0}$ is orthonormal, 
in which case the relations simplify to give 
\begin{equation}
\label{eq:relort}
\begin{split}
 r_{0,i} &= \sum_{j=1}^{n} [x_{j},[x_{j},x_{i}]] - \frac{1}{2} \sum_{a,b = 1}^{s} \Gamma^{i}_{a,b} [z_{a},z_{b}], 
 \\
 r_{1,a} &= \sum_{i=1}^{n} \sum_{b = 1}^{s} \Gamma_{a,b}^{i} [x_{i},z_{b}],  
\end{split}
\end{equation} 
which we will assume from now on. 
Therefore, if the super vector space $V$ has super dimension $(n,s)$, we may also denote the super Yang-Mills algebras by 
$\ym(n,s)^{\Gamma}$ and by $\YM(n,s)^{\Gamma}$, respectively. 
We may also write $R(n,s)^{\Gamma}$ instead of $R(V,g)^{\Gamma}$. 
Note that they generalize the Yang-Mills algebras defined by A. Connes and M. Dubois-Violette in \cite{CD02}, since $\ym(n) = \ym(n,0)^{0}$, 
and have been previously considered by M. Movshev and A. Schwarz in \cite{Mov05} and \cite{MS06}. 
We would also like to mention that, as far as we know, there is no direct relation between the previously defined super algebra and the superized versions of the Yang-Mills algebra defined in \cite{CD07}, Section 1.4, 
and, more generally, in \cite{HKL08}, Example 3.2. 
We are aware of the fact that the superized version of the Yang-Mills algebra in \cite{CD07}, Section 1.4, was also called super Yang-Mills algebra there, 
but we decided to use the same name just due to the connection to supersymmetric Yang-Mills theory in physics. 

We say that the super Yang-Mills $\ym(n,s)^{\Gamma}$ is \emph{equivariant} if $V_{1}$ is a representation of the Lie algebra $\so(V_{0},g)$ 
such that the map $\Gamma : S^{2}V^{*}_{1} \rightarrow V_{0}$ is $\so(V_{0},g)$-equivariant, and there also exists an $\so(V_{0},g)$-equivariant map 
$\tilde{\Gamma} : S^{2}V_{1} \rightarrow V_{0}$, which satisfy the following condition. 
Rewritting $\Gamma$ and $\tilde{\Gamma}$ as elements $\gamma \in \Hom (V_{0},\Hom(V_{1}^{*},V_{1}))$ and $\tilde{\gamma} \in \Hom (V_{0},\Hom(V_{1},V_{1}^{*}))$ 
defined by $z^{*}_{2}(\gamma(v)(z^{*}_{1})) = g(\Gamma(z^{*}_{1},z^{*}_{2}),v)$, for all $v \in V_{0}$ and $z_{1}^{*}, z_{2}^{*} \in V_{1}^{*}$, and 
by $\tilde{\gamma}(v)(z_{1})(z_{2}) = g(\tilde{\Gamma}(z_{1},z_{2}),v)$, for all $v \in V_{0}$ and $z_{1}, z_{2} \in V_{1}$, respectively, 
the assumption reads as follows
\begin{equation}
\label{eq:relgammas}
\begin{split}
   \tilde{\gamma}(v) \circ \gamma(v) &= g(v,v) \id_{V_{1}^{*}},     
   \\
   \gamma(v) \circ \tilde{\gamma}(v) &= g(v,v) \id_{V_{1}},      
\end{split}
\end{equation} 
for all $v \in V_{0}$.
In particular, this implies that $\tilde{\Gamma}$ is uniquely determined from $\Gamma$. 
If we denote $\tilde{\Gamma}^{i,a,b} = x_{i}^{*}(\Gamma(z_{a},z_{b}))$, for $i = 1, \dots, n$ and $a, b = 1, \dots, s$, we notice that 
the conditions \eqref{eq:relgammas} can be rewritten as 
\begin{equation}
\label{eq:relgammasind}
   \sum_{b = 1}^{s} (\Gamma^{i}_{a,b} \tilde{\Gamma}^{j,b,c} + \Gamma^{j}_{a,b} \tilde{\Gamma}^{i,b,c}) = 2 g^{i,j} \delta_{a,c}. 
\end{equation} 
We remark the easy fact that an equivariant super Yang-Mills algebra with $s \neq 0$ satisfies \textit{a fortiori} that $s(s+1)/2 \geq n$, 
for $V_{0}$ is an irreducible $\so(V_{0},g)$-module and $\Gamma \neq 0$. 

Note that if the super Yang-Mills is equivariant, then the action of $\SO(V_{0},g)$ (and so the action of $\so(V_{0},g)$) 

on $V$ induces an action by automorphisms of $\SO(V_{0},g)$ (and hence an action by derivations of $\so(V_{0},g)$) both on the tensor super algebra $TV$ 
and in the free super Lie algebra $\f(V)$, which preserves the corresponding ideal $\cl{R(n,s)^{\Gamma}}$ generated by $R(n,s)^{\Gamma}$ in the tensor algebra 
$TV$ and in the free super Lie algebra $\f(V)$, respectively. 
As a consequence, we get an action by automorphisms of $\SO(V_{0},g)$, and an action by derivations of $\so(V_{0},g)$, 
on both super Yang-Mills algebras $\YM(n,s)^{\Gamma}$ and $\ym(n,s)^{\Gamma}$. 

Even though in the examples coming from physics the super Yang-Mills algebra is equivariant, it will be useful to consider a weaker notion. 
We say that $\Gamma$, or even the super Yang-Mills algebra, is \emph{nondegenerate} 
if $n \neq 0$ and there exists a nonzero linear form $\lambda \in V_{0}^{*}$ such that 
$\lambda \circ \Gamma : V_{1}^{*} \otimes V_{1}^{*} \rightarrow k$ is nondegenerate, if $s \neq 0$. 
We remark that the super Yang-Mills algebra is nondegenerate if $s = 0$. 
Also, note that condition \eqref{eq:relgammasind} implies that an equivariant super Yang-Mills algebra with $s \neq 0$ is always nondegenerate (in fact each matrix $\Gamma^{i}$ 
is invertible). 
From now on, we shall assume that the super Yang-Mills algebras we are considering are nondegenerate. 
Without loss of generality, we shall furthermore suppose that, if $s \neq 0$, $\lambda = x_{1}^{*}|_{V_{0}}$ and that $\B_{1}^{*}$ is orthonormal with respect to 
$x_{1}^{*}|_{V_{0}} \circ \Gamma$. 
This hypothesis of nondegeneracy is though not necessary for most of this subsection and the next one, except for some minor indicated cases, 
but it will be necessary (and assumed) from Subsection \ref{subsec:anodesc} on. 

As mentioned before, the case where $s = 0$ and $n \geq 2$ has been previously studied in \cite{CD02} (see also \cite{HS10}). 
If $n=1$ and $s=0$, then the (super) Yang-Mills algebra is just the one-dimensional abelian Lie algebra, which was not studied in the mentioned articles because of its simplicity. 
Moreover, if $n = 2$ and $s=0$, $\ym(2)$ is isomorphic to the \emph{Heisenberg Lie algebra} $\h_{1}$, with generators $x,y,z$, 
and relations $[x,y] = z$, $[x,z] = [y,z] = 0$, and, for $n \geq 3$, $\ym(n)$ is an infinite dimensional Lie algebra (see \cite{HS10}, Rem. 3.14). 
On the other hand, for $n=0$ (so $s > 0$), the super Yang-Mills algebra is just the free super Lie algebra generated by the odd elements $z_{1}, \dots, z_{s}$. 
The nondegeneracy assumption on the $\Gamma$ also implies that $\ym(1,s) \simeq k.x_{1} \oplus \f(V(1))/\cl{\sum_{a=1}^{s} [z_{a},z_{a}]}$, and then 
$\YM(1,s) \simeq k[x_{1}] \otimes k\cl{z_{1},\dots,z_{s}}/\cl{\sum_{a=1}^{s} z_{a}^{2}}$. 
In particular, $\ym(1,1)^{\Gamma}$ is a supercommutative super Lie algebra of super dimension $(1,1)$, 
and $\ym(1,2)^{\Gamma}$ has super dimension $(2,2)$, with basis $x_{1}$, $z_{1}$, $z_{2}$ and $[z_{2},z_{2}]$, where $[z_{1},z_{1}]=-[z_{2},z_{2}]$, and all other brackets vanish. 
We remark that, since the case where $s=0$ is trivial, 
we will focus ourselves on indices $(n,s) \in \NN \times \NN_{0}$, unless otherwise stated.  

As noted in \cite{MS06}, it can be also useful to consider $\ym(n,s)^{\Gamma}$ as an $\NN$-graded Lie algebra, where the elements $x_{i}$ are of degree $2$, 
for all $i = 1, \dots, n$, and the elements $z_{a}$ of degree $3$, for all $a = 1, \dots, s$. 
Similarly, taking the grading induced by the previous definitions, the associative super Yang-Mills algebra $\YM(n,s)^{\Gamma}$ can also be regarded as an 
$\NN_{0}$-graded algebra. 
These gradings for both the Lie and associative versions of the super Yang-Mills algebra are exactly the \emph{special gradings} 
(opposed to the \emph{usual} ones) considered in \cite{HS10}, Section 2, and \cite{HS10h}, Section 2.1, when $s=0$. 
We remark that the underlying super Lie algebra and super algebra of these latter definitions yield the ones given at the beginning. 
Moreover, when the super Yang-Mills algebra $\ym(n,s)^{\Gamma}$ is seen as a graded Lie algebra, we may consider the descending sequence of graded ideals $\{ F^{\bullet}\ym(n,s)^{\Gamma} \}_{\bullet \in \NN_{0}}$, 
where $F^{j}\ym(n,s)^{\Gamma}$ is the graded vector subspace of $\ym(n,s)^{\Gamma}$ given by elements of degree greater than or equal to $j + 2$. 
It is the quotient under the canonical projection of the descending filtration $F^{\bullet}\f(V(n,s))$ of the graded free Lie algebra 
given by elements of degree greater than or equal to $j + 2$. 
This also induces a descending sequence of ideals of the underlying super Lie algebra of either the graded free Lie algebra $\f(V(n,s))$ 
or the graded Lie algebra $\ym(n,s)^{\Gamma}$ considered above. 

\begin{remark}
\label{rem:phys}
We remark that this super algebra appears naturally when studying supersymmetric gauge field theories. 
We will only recall what we need for our explanation, and we refer to \cite{DF99} and \cite{DF99s} for a complete account on this subject. 
Let $\check{M}$ denote an $n$-dimensional \emph{Minkowski space} with real vector space of translations $V$ of dimension $n$ 
and metric $g$ determining a cone $C$ of time-like vectors in $V$, $S$ a real spinorial representation of dimension $s$ of the spin group $\mathrm{Spin}(V,g)$ 
and $\Gamma : S^{*} \otimes S^{*} \rightarrow V$ a symmetric morphism of representations of $\mathrm{Spin}(V,g)$, which is \emph{positive definite}, 
\textit{i.e.} $\Gamma(s^{*},s^{*}) \in \bar{C}$, for all $s^{*} \in S^{*}$, and $\Gamma(s^{*},s^{*}) = 0$ only if $s^{*} = 0$. 
Since any complex vector bundle $E$ of rank $m$ over $\check{M}$ is trivial, for $\check{M}$ is contractible, 
every connection on such bundle is given by a $\gl_{m}(\CC)$-valued $1$-form $\sum_{i=1}^{n} A_{i} dx^{i}$, 
and the corresponding covariant derivative is $\nabla_{i} = \partial_{i} + A_{i}$, for $i = 1, \dots, n$. 
We also recall that a \emph{dual spinor field} $\lambda$ with values in the Lie algebra $\gl_{m}(\CC)$ is a morphism from $\check{M}$ to 
$S^{*} \otimes \gl_{m}(\CC)$, so it can be decomposed in components $\lambda_{a}$ from $\check{M}$ to $\gl_{m}(\CC)$, for $a = 1, \dots, s$. 
We may thus see both set of fields $\nabla_{i}$ and $\lambda_{a}$ as sections of the endomorphism bundle of the complex super vector bundle 
$E \otimes_{\RR} \Lambda^{\bullet}_{\RR} S^{*}$ on $\check{M}$, where $S$ is considered to be in degree $1$, 
$\nabla_{i}$ is even, for all $i = 1, \dots, n$, and $\lambda^{a}$ is odd, for all $a = 1, \dots, s$. 
The associated set of \emph{super Yang-Mills equations} is the supersymmetric extension of the usual Yang-Mills equations 
and it is given by 
\begin{equation}
\label{eq:ym}
\begin{split}
 \sum_{j,l,m=1}^{n} g^{i,j} g^{l,m}[\nabla_{j},[\nabla_{l},\nabla_{m}]] &= \frac{1}{2} \sum_{a,b = 1}^{s} \Gamma^{i}_{a,b} [\lambda_{a},\lambda_{b}], 
 \\
 \sum_{i=1}^{n} \sum_{b = 1}^{s} \Gamma_{a,b}^{i} [\nabla_{i},\lambda_{b}] &= 0,  
\end{split}
\end{equation}
for $i = 1, \dots, n$ and $a = 1, \dots, s$, respectively (\textit{cf.} \cite{DF99s}, (6.23)). 
We remark that the previous identities are considered in the super Lie algebra of endomorphisms of sections of the complex super vector bundle 
$E \otimes_{\RR} \Lambda^{\bullet}_{\RR} S^{*}$, so they give a representation of the equivariant super Yang-Mills algebra $\ym(n,s)^{\Gamma}$. 
\end{remark}

The Lie ideal $\tym(n,s)^{\Gamma} = F^{1}\ym(n,s)^{\Gamma} = [\ym(n,s)^{\Gamma},\ym(n,s)^{\Gamma}] + \ym(n,s)^{\Gamma}_{1}$ 
of the super Lie algebra $\ym(n,s)^{\Gamma}$ will be also important in the sequel. 
It is obvious to see that it is also a graded Lie ideal of $\ym(n,s)^{\Gamma}$, when considered as a graded Lie algebra. 
Notice that $\tym(n,0)^{0}$ coincides with the ideal $\tym(n)$ of $\ym(n) = \ym(n,0)^{0}$ considered in \cite{HS10}. 
We remark that $\ym(n,s)^{\Gamma}/\tym(n,s)^{\Gamma} \simeq V(n)_{0}$ is the abelian (super) Lie algebra of super dimension $(n,0)$. 
Moreover, we shall deal with the universal enveloping algebra of the Lie ideal $\tym(n,s)^{\Gamma}$, which will be denoted $\TYM(n,s)^{\Gamma}$, 
which can be regarded either as a super algebra or as a graded algebra. 
We shall also consider the bigger Lie ideal $\hat{\tym}(n,s)^{\Gamma} = \tym(n,s)^{\Gamma} \oplus \bigoplus_{i=3}^{n} k.x_{i}$, 
which satisfies that $\ym(n,s)^{\Gamma}/\hat{\tym}(n,s)^{\Gamma} \simeq V(2)_{0}$ is supercommutative, and its enveloping algebra 
$\hat{\TYM}(n,s)^{\Gamma} = \U(\hat{\tym}(n,s)^{\Gamma})$. 
Occasionally, we will omit the indices $(n,s)$ and $\Gamma$ for the previously defined spaces (and also for the defined below) 
in order to simplify the notation if it is clear from the context.

We would like to make a further comment on the relationship between different super Yang-Mills algebras. 
One first notes that the canonical projection $V(n,s) \rightarrow V(n)_{0}$ of graded vector spaces induces a surjective morphism of graded algebras 
$TV(n,s) \rightarrow TV(n)_{0}$, which sends the even super Yang-Mills relations $r_{0,i}$ to the usual Yang-Mills relations $r_{i}$ described in \cite{HS10}, 
and the odd super Yang-Mills relations $r_{1,a}$ to zero. 
Hence, we obtain a surjective morphism of graded algebras $\YM(n,s)^{\Gamma} \rightarrow \YM(n)$, 
where we recall that the Yang-Mills algebras are provided with the special grading.  
In an analogous manner, we have a surjective morphism of graded Lie algebras $\ym(n,s)^{\Gamma} \rightarrow \ym(n)$, which obviously maps 
$\tym(n,s)^{\Gamma}$ onto $\tym(n)$. 
We would like to point out however, that we do not know of any \textit{a priori} given morphism between different super Yang-Mills algebras with $s \neq 0$ 
when the corresponding $\Gamma$'s are arbitrary. 

\subsection{A superpotential formulation}
\label{subsec:pot}

We remark that relations \eqref{eq:relort} can be obtained from the superpotential element given by the class of the homogeneous element 
(resp., of degree $8$ if $\ym(n,s)^{\Gamma}$ is considered as a graded Lie algebra)
\begin{equation}
\label{eq:sup}
     W = - \frac{1}{4} \sum_{i,j = 1}^{n} [x_{i},x_{j}][x_{i},x_{j}] + \frac{1}{2} \sum_{i = 1}^{n} \sum_{a,b = 1}^{s} \Gamma_{a,b}^{i} z_{a}[x_{i},z_{b}]     
\end{equation} 
in $HH_{0}(TV) = TV/[TV,TV]$, where we recall that we are considering the super (resp. graded) commutator space $[TV,TV]$ of the tensor algebra, \textit{i.e.} 
the super (resp., graded) vector space spanned by the elements $u v - (-1)^{|u||v|} v u$, for homogeneous elements $u, v \in TV$. 
In this case, the \emph{cyclic derivative with respect to a generator} $v \in \B$, where $\B$ is a basis of (homogeneous elements of) $V$, is the map from $TV/[TV,TV]$ to $TV$ given by the 
$k$-linear extension of the following: for an element $w = v_{1} \dots v_{r}$, where each $v_{i} \in \B$ for all $i = 1, \dots, r$, we consider 
\[     \frac{\partial \bar{w}}{\partial v} = \sum_{i : v_{i} = v} (-1)^{(|v_{i}| + \dots + |v_{n}|) (|v_{1}| + \dots + |v_{i-1}|)}
                                              v_{i+1} \dots v_{n} v_{1} \dots v_{i-1}.     \]
It can be easily checked that in our case $r_{0,i}$ is the cyclic derivative of $W$ with respect to $x_{i}$, and $r_{1,a}$ is the cyclic derivative of $W$ with respect to $z_{a}$. 

We have the following proposition, which we suppose it should be well-known, since it is only the algebraic analogous of the physical folkloric result that states that  
the Euler-Lagrange equations associated to an action are invariant under a collection of (super)symmetries if the action is invariant under such (super)symmetries.
\begin{proposition}
\label{prop:simpot}
Let $V$ be a graded vector space $V$, and let $\g$ be a graded Lie algebra acting by homogeneous derivations 
on the corresponding tensor algebra $TV$.  
Regarding $TV$ as a graded Lie algebra, we immediately note that the action of $\g$ on $TV$ induces an action on (the abelianization) $HH_{0}(TV) = TV/[TV,TV]$. 
If a homogeneous element $\bar{W} \in HH_{0}(TV)$ belongs to the invariant space $HH_{0}(TV)^{\g}$, then the two-sided ideal generated by the 
cyclic derivatives of $W$ (with respect to any basis of $V$) is preserved by $\g$. 
\end{proposition}
\noindent\textbf{Proof.}
It suffices to prove the statement for a unique homogeneous derivation $d$ of $TV$ of degree $|d|$, and call $\bar{d}$ the induced map on $HH_{0}(TV)$. 
Let $\B$ be a (homogeneous) basis of $V$. 
Fix an element $v$ of the basis $\B$ of $V$, and for each $v' \in \B$, 
consider the monomial elements $r_{l,v'}^{m}, s_{l,v'}^{m} \in TV$ defined by $d(v') = \sum_{m \in M_{v'}} r_{l,v'}^{m} v s_{l,v'}^{m}$, for all $l \in L_{v',m}$, and 
\[     \frac{\partial (\bar{d}(\bar{v}'))}{\partial v} = \sum_{m \in M_{v'}} \sum_{l \in L_{v',m}} (-1)^{|r_{l,v'}^{m}|(|v|+|s_{l,v'}^{m}|)} 
       s_{l,v'}^{m} r_{l,v'}^{m}     \]
hold, where $M_{v'}$ and $L_{v',m}$ are sets of indices. 
Set $\B' = \B \setminus \{ v \}$. 
We would like to explain the notation, which could seem cumbersome at first glance. 
The element $d(v')$ is a sum of monomials $M_{m}$ in $TV$, indexed by $m \in M_{v'}$, and each of these monomials can be written as $M_{m} = r_{l,v'}^{m} v s_{l,v'}^{m}$, for $l \in L_{v',m}$, 
which are all the different ways of writting $M_{m}$ such that the formula for the cyclic derivative can be applied to give 
\[     \frac{\partial (\bar{d}(\bar{M_{m}}))}{\partial v} = \sum_{l \in L_{v',m}} (-1)^{|r_{l,v'}^{m}|(|v|+|s_{l,v'}^{m}|)} s_{l,v'}^{m} r_{l,v'}^{m}.     \] 

We claim that the following identity holds 
\begin{equation}
\label{eq:difi}
      \frac{\partial (\bar{d}(\bar{W}))}{\partial v} = (-1)^{|d||v|} d\Big(\frac{\partial \bar{W}}{\partial v}\Big) 
       + \sum_{v' \in \B} \sum_{m \in M_{v'}} \sum_{l \in L_{v',m}} (-1)^{|r_{l,v'}^{m}|(|W|+|s_{l,v'}^{m}| + |v|-|v'|)} s_{l,v'}^{m} 
         \frac{\partial \bar{W}}{\partial v'} r_{l,v'}^{m}.     
\end{equation}

First, note that it only suffices to prove the previous equality for $W$ a monomial, 
since if we write $W = \sum_{j \in J} W_{j}$, where $W_{j}$ are monomials (of the same degree as $W$) in $TV$, the sum of each of the identities \eqref{eq:difi} for 
$W_{j}$ gives \eqref{eq:difi} for $W$. 
Then, we assume that $W$ is a monomial of the form $\prod_{i = 1}^{N} v_{i}$, for $N \in \NN$, such that 
$v_{i} \in \B$. 

Define $d_{i,j}$ to be $\sum_{l=i}^{j} |v_{l}|$, for $1 \leq i \leq j \leq N$, and zero else.

On the one hand, we have that 
\[     \frac{\partial \bar{W}}{\partial v} = \sum_{i : v_{i} = v} (-1)^{d_{1,i-1}d_{i,N}} (\prod_{j=i+1}^{N} v_{j}) (\prod_{j=1}^{i-1} v_{j}),     \]
which yields 
\begin{equation}  
\label{eq:dpW}
\begin{split}
       d\Big(\frac{\partial \bar{W}}{\partial v}\Big) 
       =& \sum_{i : v_{i} = v} (-1)^{d_{1,i-1}d_{i,N}} \Big(\sum_{p = i+1}^{N} (-1)^{|d|d_{i+1,p-1}} 
          (\prod_{j=i+1}^{p-1} v_{j}) d(v_{p}) (\prod_{j=p+1}^{N} v_{j}) (\prod_{j=1}^{i-1} v_{j}) 
       \\
       &+ \sum_{p = 1}^{i-1} (-1)^{|d|(d_{i+1,N}+d_{1,p-1})} 
          (\prod_{j=i+1}^{N} v_{j}) (\prod_{j=1}^{p-1} v_{j}) d(v_{p}) (\prod_{j=p+1}^{i-1} v_{j}) \Big).
\end{split}
\end{equation}

On the other hand, from 
\begin{align*}
       d(W) &= \sum_{p=1}^{N} (-1)^{|d|d_{1,p-1}} (\prod_{j=1}^{p-1} v_{j}) d(v_{p}) (\prod_{j=p+1}^{N} v_{j})   
       \\
            &= \sum_{p=1}^{N} \sum_{m \in M_{v_{p}}} (-1)^{|d|d_{1,p-1}} (\prod_{j=1}^{p-1} v_{j}) r_{l,v_{p}}^{m} v s_{l,v_{p}}^{m} (\prod_{j=p+1}^{N} v_{j}), 
\end{align*}
we get that
\begin{align*}
   \frac{\partial (\bar{d}(\bar{W}))}{\partial v} 
   = \sum_{p=1}^{N} (-1)^{|d|d_{1,p-1}} 
      &\Big(\sum_{i < p : v_{i} = v} (-1)^{d_{1,i-1}(|d|+d_{i,N})} (\prod_{j=i+1}^{p-1} v_{j}) d(v_{p}) (\prod_{j=p+1}^{N} v_{j}) (\prod_{j=1}^{i-1} v_{j}) 
    \\
      &+ \sum_{i > p : v_{i} = v} (-1)^{d_{i,N}(|d|+d_{1,i-1})} (\prod_{j=i+1}^{N} v_{j}) (\prod_{j=1}^{p-1} v_{j}) d(v_{p}) (\prod_{j=p+1}^{i-1} v_{j}) 
    \\
      &+ \sum_{m \in M_{v_{p}}} \sum_{l \in L_{v_{p}}} (-1)^{(d_{1,p-1}+|r_{l,v_{p}}^{m}|)(|v|+|s_{l,v_{p}}^{m}|+d_{p+1,N})} 
                 s_{l,v_{p}}^{m} (\prod_{j=p+1}^{N} v_{j}) (\prod_{j=1}^{p-1} v_{j}) r_{l,v_{p}}^{m}\Big).   
\end{align*} 
It is obvious to see that the first two terms of the second member of the previous identity coincide with $(-1)^{|v||d|}$ times the two terms of the last member of 
\eqref{eq:dpW}. 
It is also clear that the third term of the equation coincides with the second term of the second member of \eqref{eq:difi}, which proves the claim. 
The proposition now follows directly from the identity \eqref{eq:difi}, since $\bar{d}(\bar{W}) = 0$. 
\qed

\begin{remark}
\label{rem:simpot}
From the proof of the proposition we further see that the ideal $I \subset TV$ generated by the cyclic derivatives of a homogeneous superpotential 
$\bar{W} \in HH_{0}(TV)$ is preserved by a homogeneous derivation $d$ if and only if the cyclic derivatives of $\bar{d}(\bar{W})$ belong to $I$. 
\end{remark}

For the rest of this subsection we assume that the Yang-Mills algebra $\ym(n,s)^{\Gamma}$ is equivariant and that $V_{1}^{*}$ is an irreducible spin representation of $\so(V_{0},g)$ 
(\textit{i.e.}, a so-called \emph{minimal supersymmetry} in physical theories). 
We have the following direct consequence of the proposition. 
Let us consider the collection $d_{c}^{1}$ ($c = 1, \dots, s$) of homogeneous derivations of $TV$ of degree $1$ induced by 
\begin{align*}
   d_{c}^{1}(x_{i}) &= \sum_{d=1}^{s} \Gamma^{i}_{c,d} z_{d}, 
   \\
   d_{c}^{1}(z_{b}) &= \frac{1}{2} \sum_{d=1}^{s} \sum_{i,j=1}^{n} \tilde{\Gamma}^{i,b,d} \Gamma^{j}_{d,c} [x_{i},x_{j}]. 
\end{align*}
As expected, they are analogous to the supersymmetry transformations that have been considered in supersymmetric gauge theories 
(\textit{cf.} \cite{DF99s}, (6.8)). 
The following computation is parallel to one already known to physicists long time ago, but we provide it just for completeness (\textit{cf.} \cite{DF99s}, Thm. 6.4). 
Applying the derivation $d_{a}^{1}$ to the superpotential $\bar{W}$ given in \eqref{eq:sup}, we see that, on the one hand 
\[     d_{c}^{1}\Big(-\frac{1}{4} \sum_{i,j = 1}^{n} [x_{i},x_{j}][x_{i},x_{j}]\Big) = - \sum_{b=1}^{s} \sum_{i,j = 1}^{n} \Gamma^{j}_{c,b} [x_{i},x_{j}][x_{i},z_{b}].     \]
On the other hand, we have 
\begin{align*}
   d_{c}^{1}\Big(\frac{1}{2} \sum_{i = 1}^{n} \sum_{a,b = 1}^{s} \Gamma_{a,b}^{i} z_{a}[x_{i},z_{b}]\Big) 
   &= \frac{1}{4} \sum_{i,j,l = 1}^{n} \sum_{a,b,d = 1}^{s} 
                  \Gamma_{a,b}^{i} \tilde{\Gamma}^{j,a,d} \Gamma^{l}_{d,c} [x_{j},x_{l}] [x_{i},z_{b}] 
     - \frac{1}{2} \sum_{i = 1}^{n} \sum_{a,b,d = 1}^{s} \Gamma_{a,b}^{i} \Gamma^{i}_{c,d} z_{a} [z_{d},z_{b}] 
   \\
   &- \frac{1}{4} \sum_{i,j,l = 1}^{n} \sum_{a,b,d = 1}^{s} \Gamma_{a,b}^{i} \tilde{\Gamma}^{j,b,d} \Gamma^{l}_{d,c} 
           z_{a}[x_{i}, [x_{j},x_{l}]]
   \\
   &= - \frac{1}{2} \sum_{i,j,l = 1}^{n} \sum_{a,b,d = 1}^{s}
                  \Gamma_{a,b}^{i} \tilde{\Gamma}^{j,b,d} \Gamma^{l}_{d,c} [x_{i},[x_{j},x_{l}]] z_{a}
     - \frac{1}{2} \sum_{i = 1}^{n} \sum_{a,b,d = 1}^{s} \Gamma_{a,b}^{i} \Gamma^{i}_{c,d} z_{a} [z_{d},z_{b}] 
   \\
   &+ \frac{1}{4} \sum_{i,j,l = 1}^{n} \sum_{a,b,d = 1}^{s} \Gamma_{a,b}^{i} \tilde{\Gamma}^{j,b,d} \Gamma^{l}_{d,c} 
           [x_{i}, [x_{j},x_{l}]z_{a}],
\end{align*}
where we have used that $[x_{j},x_{l}][x_{i},z_{a}] = [x_{i},[x_{j},x_{l}]z_{a}] - [x_{i},[x_{j},x_{l}]]z_{a}$ in the first term of the third member. 

We now consider the following collection of homogeneous element of $TV$ (for $a,c = 1, \dots, s$)
\begin{align*}
      X_{a,c} &= \frac{1}{2} \sum_{i,j,l = 1}^{n} \sum_{b,d = 1}^{s}
                  \Gamma_{a,b}^{i} \tilde{\Gamma}^{j,b,d} \Gamma^{l}_{d,c} [x_{i},[x_{j},x_{l}]] 
               = \frac{1}{2} \sum_{i,j,l = 1}^{n} \sum_{b,d = 1}^{s}
                  \Gamma_{a,b}^{i} \tilde{\Gamma}^{j,b,d} \Gamma^{l}_{d,c} ([x_{l},[x_{j},x_{i}]] + [x_{j},[x_{i},x_{l}]]) 
      \\
               &=  \frac{1}{2} \sum_{i,j,l = 1}^{n} \sum_{b,d = 1}^{s}
                  (\Gamma_{c,b}^{i} \tilde{\Gamma}^{j,b,d} \Gamma^{l}_{d,a} + \Gamma_{a,b}^{j} \tilde{\Gamma}^{i,b,d} \Gamma^{l}_{d,c}) [x_{i},[x_{j},x_{l}]] 
      \\       
               &=  X_{c,a} + \frac{1}{2} \sum_{i,j,l = 1}^{n} \sum_{d = 1}^{s}
                  (- \sum_{b = 1}^{s} \Gamma_{a,b}^{i} \tilde{\Gamma}^{j,b,d} \Gamma^{l}_{d,c} + 2 \Gamma_{d,c}^{l} \delta_{i,j} \delta_{a,d}) [x_{i},[x_{j},x_{l}]] 
      \\
               &= X_{c,a} - X_{a,c} + S_{a,c},      
\end{align*}
where we have used \eqref{eq:relgammasind} in the penultimate equality and we define
\[     S_{a,c} = \sum_{i,l=1}^{n} \Gamma_{a,c}^{l} [x_{i},[x_{i},x_{l}]].     \]
This implies that $2 X_{a,c} = X_{c,a} + S_{a,c}$. 
Using that $S_{a,c}$ is symmetric for the interchange of indices, a trivial computation implies that $X_{a,c}$ is also so, 
which in turn implies that $X_{a,c} = S_{a,c}$. 
As a consequence, we have that 
\begin{align*}
   d_{c}^{1}\Big(\frac{1}{2} \sum_{i = 1}^{n} \sum_{a,b = 1}^{s} \Gamma_{a,b}^{i} z_{a}[x_{i},z_{b}]\Big) 
   &= - \sum_{i,l=1}^{n} \sum_{a=1}^{s} \Gamma_{a,c}^{l} [x_{i},[x_{i},x_{l}]] z_{a}
     - \frac{1}{2} \sum_{i = 1}^{n} \sum_{a,b,d = 1}^{s} \Gamma_{a,b}^{i} \Gamma^{i}_{c,d} z_{a} [z_{d},z_{b}] 
   \\
   &+ \frac{1}{4} \sum_{i,j,l = 1}^{n} \sum_{a,b,d = 1}^{s} \Gamma_{a,b}^{i} \tilde{\Gamma}^{j,b,d} \Gamma^{l}_{d,c} [x_{i}, [x_{j},x_{l}]z_{a}]
   \\
   &= \sum_{i,l=1}^{n} \sum_{a=1}^{s} \Gamma_{a,c}^{l} [x_{i},x_{l}] [x_{i},z_{a}]
     - \frac{1}{2} \sum_{i = 1}^{n} \sum_{a,b,d = 1}^{s} \Gamma_{a,b}^{i} \Gamma^{i}_{c,d} z_{a} [z_{d},z_{b}] 
   \\
   &+ \frac{1}{4} \sum_{i,j,l = 1}^{n} \sum_{a,b,d = 1}^{s} \Gamma_{a,b}^{i} \tilde{\Gamma}^{j,b,d} \Gamma^{l}_{d,c} [x_{i}, [x_{j},x_{l}]z_{a}] 
     - \sum_{i,l=1}^{n} \sum_{a=1}^{s} \Gamma_{a,c}^{l} [x_{i},[x_{i},x_{l}] z_{a}].
\end{align*}
The last two terms belong to $HH_{0}(TV)$, so we may discard them when considering $\bar{d}_{c}^{1}(\bar{W})$. 
In consequence, summing up, we obtain that $\bar{d}_{c}^{1}(\bar{W})$ is the class of the element 
\[     - \frac{1}{2} \sum_{i = 1}^{n} \sum_{a,b,d = 1}^{s} \Gamma_{a,b}^{i} \Gamma^{i}_{c,d} z_{a} [z_{d},z_{b}],     \] 
whose cyclic derivatives belong to the ideal generated by the cyclic derivatives of $W$ if and only if it vanishes. 
Therefore, using the Jacobi identity and Proposition \ref{prop:simpot} (or, more precisely, Remark \ref{rem:simpot}), 
the set of homogeneous derivations $\{ d_{c}^{1} \}_{c=1,\dots,s}$ preserve the ideal $\cl{R(n,s)^{\Gamma}}$ if and only if 
the quartic form $\sum_{i=1}^{n} (\Gamma_{a,b}^{i} \Gamma_{c,d}^{i} + \Gamma_{a,c}^{i} \Gamma_{b,d}^{i} + \Gamma_{a,d}^{i} \Gamma_{b,c}^{i})$ 
(\textit{i.e.} the one given by $z^{*} \mapsto g(\Gamma(z^{*},z^{*}),\Gamma(z^{*},z^{*}))$, for $z^{*} \in V_{1}^{*}$) vanishes. 
As explained in \cite{DF99s}, \S 6.1, this happens for $n = 3, 4, 6, 10$. 
This also gives a simpler proof of \cite{Mov05}, Prop. 20 (besides the other implication which was not stated there). 

\subsection{Another description of the super Yang-Mills algebra}
\label{subsec:anodesc}

Define $\h(n,s)$ to be super Lie algebra generated by the super vector space $U(n,s) = U(n,s)_{0} \oplus U(n,s)_{1}$, where 
$U(n,s)_{0} = \mathrm{span}_{k} \cl{q_{2}, \dots, q_{n}, p_{2}, \dots, p_{n}}$ and $U(n,s)_{1} = \mathrm{span}_{k} \cl{z_{1}', \dots, z_{s}'}$, 
with the relation space given  by $R(n,s)' = k.(\sum_{i=2}^{n} [q_{i},p_{i}] + \frac{1}{2} \sum_{a=1}^{s} [z_{a}',z_{a}'])$, 
\textit{i.e.} $\h(n,s) = \f(U(n,s))/\cl{R(n,s)'}$. 
The (super) Lie algebra $\a = k.d$, with $d$ even, acts by (even) derivations on $\h(n,s)$ as follows: 
\begin{align*}
   d(q_{i}) &= p_{i}, 
   \\
   d(p_{i}) &= - \sum_{j=2}^{n} [q_{j},[q_{j},q_{i}]] + \frac{1}{2} \sum_{a,b = 1}^{s} \Gamma_{a,b}^{i} [z_{a}',z_{b}'], 
   \\
   d(z_{a}') &= - \sum_{j=2}^{n} \sum_{b=1}^{s} \Gamma_{a,b}^{j} [q_{j},z_{b}'],
\end{align*}
for all $i = 2, \dots, n$ and $a = 1, \dots, s$. 
Note that, for the previous action to be well-defined, $\Gamma$ need not be nondegenerate. 
We now easily obtain a morphism of super Lie algebras 
\[     \psi : \ym(n,s)^{\Gamma} \rightarrow \a \ltimes \h(n,s)     \]
given by 
\begin{align*}
   x_{1} &\mapsto d,
   \\
   x_{i} &\mapsto q_{i},
   \\
   z_{a} &\mapsto z_{a}',
\end{align*}
for all $i = 2, \dots, n$ and $a = 1, \dots, s$. 
The morphism is well-defined because of the assumption that $\B_{1}^{*}$ is orthonormal with respect to $x_{1}^{*}|_{V_{0}} \circ \Gamma$. 
Moreover, it is bijective with inverse given by the morphism of super Lie algebras 
\[     \psi^{-1} : \a \ltimes \h(n,s) \rightarrow \ym(n,s)^{\Gamma},     \]
defined as 
\begin{align*}
   d &\mapsto x_{1},
   \\
   q_{i} &\mapsto x_{i},
   \\
   p_{i} &\mapsto [x_{1},x_{i}],
   \\
   z_{a}' &\mapsto z_{a},  
\end{align*}
for all $i = 2, \dots, n$ and $a = 1, \dots, s$. 
Therefore, we have that, under the assumption that $\Gamma$ is nondegenerate, $\ym(n,s)^{\Gamma} \simeq \a \ltimes \h(n,s)$, 
which further yields that $\YM(n,s)^{\Gamma} \simeq \U(\a) \# H(n,s)$, where $H(n,s) = \U(\h(n,s))$ (\textit{cf.}~\cite{MS06}, Prop. 11). 
In particular, we may regard $\h(n,s)$ canonically included in $\ym(n,s)^{\Gamma}$ (using the morphism $\psi^{-1}$). 
We would also like to point out that all previous isomorphisms also hold if we consider $\ym(n,s)^{\Gamma}$ as a graded Lie algebra with generators $x_{i}$ in degree $2$ and generators $z_{a}$ in degree $3$, 
$\a$ as the free graded Lie algebra with generator $d$ of degree $2$, and $\h(n,s)$ as a graded Lie algebra with generators $q_{i}$ of degree $2$, 
$p_{i}$ of degree $4$ and $z_{a}'$ of degree $3$, with the given relation (of degree $6$).  

\section{Several homological computations}
\label{sec:hom}

\subsection{A Koszul-like projective resolution}
\label{subsec:projres}

In this subsection we shall provide the minimal projective resolution of the left $\YM(n,s)^{\Gamma}$-module $k$, 
when $\YM(n,s)^{\Gamma}$ is considered as a graded algebra, where we suppose that $(n,s) \in (\NN \times \NN_{0}) \setminus \{ (1,0), (1,1) \}$. 
We have excluded the case $(n,s) = (1,0)$, because in this case the super Yang-Mills algebra is a polynomial algebra in one (even) variable. 
Using ideas analogous to the case of Koszul algebras, we shall construct from the previous resolution the minimal projective resolution of the $\YM(n,s)^{\Gamma}$-bimodule $\YM(n,s)^{\Gamma}$, 
and derive several consequences from the explicit shape of the previous resolutions. 
Although some of the results of the first part of this subsection are mentioned in \cite{MS06}, our aim here is to give detailed proofs of the results that we will need later. 

Regarding $\YM(n,s)^{\Gamma}$ as a nonnegatively graded connected algebra with the grading stated in Subsection \ref{subsec:defi}, 
and since $k$ is a bounded below graded module over $\YM(n,s)^{\Gamma}$, 
we know that there exists a minimal projective resolution of $k$ over $\YM(n,s)^{\Gamma}$, which is of the form 
$K_{\bullet}' = \YM(n,s)^{\Gamma} \otimes_{k} \Tor_{\bullet}^{\YM(n,s)^{\Gamma}}(k,k)$ 
(see \cite{Ber08}, Th\'eo. 1.11, Prop. 2.3). 
In order to obtain this minimal projective resolution, we shall proceed as follows. 
First, we will compute the homology groups $H_{\bullet}(\h(n,s),k)$ together with their action of $\a$. 
From the Hochschild-Serre spectral sequence $E_{p,q}^{2} = H_{p}(\a,H_{q}(\h(n,s),k)) \Rightarrow H_{p+q}(\ym(n,s)^{\Gamma},k)$, we shall obtain 
the homology groups $H_{\bullet}(\ym(n,s)^{\Gamma},k)$. 
These computations will in turn give the minimal projective resolution of $k$ over $\YM(n,s)^{\Gamma}$. 
We remark that all algebras $\a$, $\h(n,s)$ and $\ym(n,s)^{\Gamma}$ in these computations are considered as graded Lie algebras with the grading explained at the end of Subsections \ref{subsec:anodesc}. 

We first compute the homology groups $H_{\bullet}(\h(n,s),k)$ together with their action of $\a$. 
The minimal projective resolution of $k$ as a module over $H(n,s) = \U(\h(n,s))$ is of the form 
\begin{equation}
\label{eq:hns}
     0 \rightarrow H(n,s) \otimes R(n,s)' \overset{d_{2}'}{\rightarrow} H(n,s) \otimes U(n,s) \overset{d_{1}'}{\rightarrow} 
                     H(n,s) \overset{\epsilon_{\h(n,s)}}{\rightarrow} k \rightarrow 0,    
\end{equation} 
where we recall that $d_{2}'$ is given by the restriction of the map $H(n,s) \otimes U(n,s)^{\otimes 2} \rightarrow H(n,s) \otimes U(n,s)$ 
of the form $z \otimes (v \otimes w) \mapsto z.v \otimes w$, for $z \in H(n,s)$ and $v, w \in U(n,s)$, and 
$d_{1}'$ is the restriction of the multiplication of $H(n,s)$. 
Indeed, if $H(n,s)$ is considered to be generated in degree one, \textit{i.e.} $U(n,s)$ is seen to be concentrated in degree $1$, $H(n,s)$ becomes a quadratic algebra with only one relation, therefore Koszul 
(see \cite{Fr75}, Thm. 3). 
moreover, the minimal projective resolution of the $H(n,s)$-module $k$ is given by \eqref{eq:hns}, because the map $d_{2}'$ is injective. 
To prove this last statement we may note that, if $n \geq 2$ then $\U(\ym(n,s)^{\Gamma}) \supset \U(\ym(n))$ is a free extension of algebras with $\U(\ym(n))$ a
domain, and if $n = 1$ and $s \geq 2$, it follows from a simple computation using an explicit basis of the algebra $\U(\h(1,s))$ given by the Diamond Lemma (see \cite{Berg78}). 
Since the resolution \eqref{eq:hns} also holds for the special grading of $H(n,s)$, \textit{i.e.} when $\deg(q_{i}) = 2$, $\deg(p_{i}) = 4$, and $\deg(z_{a}') = 3$, for all $i = 2, \dots, n$ and $a = 1, \dots, s$, 
the claim follows. 
In particular we have canonical isomorphisms $H_{0}(\h(n,s),k) \simeq k$, $H_{1}(\h(n,s),k) \simeq U(n,s)$, $H_{2}(\h(n,s),k) \simeq R(n,s)'$ 
and the other homology groups vanish. 

In order to obtain the action of $\a$ on the homology groups we compare the previous resolution with the Chevalley-Eilenberg resolution of $k$ over 
the graded Lie algebra $\ym(n,s)^{\Gamma}$, seen as a resolution of $k$ as $\h(n,s)$-modules. 
We recall that, for a graded (or super) Lie algebra $\g$, the Chevalley-Eilenberg resolution of the trivial $\g$-module $k$ is of the form $(\U(\g) \otimes \Lambda^{\bullet} \g, d^{\mathrm{CE}}_{\bullet})$, 
where we recall that the exterior tensor products are taken in the graded (or super) sense, and the differential is given by (\textit{cf.} \cite{Ta95}, Section 1) 
\begin{multline*}
     d_{n}^{\mathrm{CE}}(z \otimes y_{1} \wedge \dots \wedge y_{n}) = \sum_{i=1}^{n} (-1)^{|y_{i}|(|y_{1}|+\dots+|y_{i-1}|) + (i-1)} z y_{i} \otimes y_{1} \wedge \dots \wedge \hat{y}_{i} \wedge \dots \wedge y_{n} 
     \\
                                                          + \sum_{1 \leq i < j \leq n} (-1)^{|y_{i}|(|y_{1}|+\dots+|y_{i-1}|) + |y_{j}|(|y_{1}|+\dots+|y_{j-1}|) + |y_{i}||y_{j}| + (i+j)} 
                                                                  z \otimes [y_{i},y_{j}] \wedge y_{1} \wedge \dots \wedge \hat{y}_{i} \wedge \dots \wedge \hat{y}_{j} \wedge \dots \wedge y_{n},   
\end{multline*}
for $y_{1}, \dots y_{n} \in \g$ homogeneous elements and $z \in \U(\g)$. 

We have the following comparison morphism of resolutions of projective $H(n,s)$-modules of $k$: 
\[
\xymatrix
{
\dots
\ar[r]
&
0 
\ar[r] 
\ar[d]
&
H \otimes R' 
\ar[r]^{d_{2}'}
\ar[d]^{\mathrm{can}}
&
H \otimes U
\ar[r]^{d_{1}'}
\ar[d]^{\mathrm{inc}}
&
H
\ar[r]^{\epsilon_{\h}}
\ar[d]^{\mathrm{inc}}
&
k
\ar[r]
\ar@{=}[d]
&
0
\\
\dots
\ar[r]^-{d_{4}^{CE}}
&
\YM \otimes \Lambda^{3} \ym
\ar[r]^-{d_{3}^{CE}}
&
\YM \otimes \Lambda^{2} \ym
\ar[r]^-{d_{2}^{CE}}
&
\YM \otimes \ym
\ar[r]^-{d_{1}^{CE}}
&
\YM
\ar[r]^{\epsilon_{\ym}}
&
k
\ar[r]
&
0
}    
\] 
where we have omitted the indices $(n,s)$ and $\Gamma$ for simplicity, 
\[     \mathrm{can}(z \otimes (\sum_{i=2}^{n} [q_{i},p_{i}] + \frac{1}{2} \sum_{a=1}^{s} [z_{a}',z_{a}'])) 
       = z \otimes (\sum_{i=2}^{n} x_{i} \wedge [x_{1},x_{i}] + \frac{1}{2} \sum_{a=1}^{s} z_{a} \wedge z_{a}),     \]
and we remark that the exterior power and the wedge product here are in the supersymmetric sense. 

Applying the functor $(\place)_{\h(n,s)}$ to both resolutions, we see that, under the previous comparison morphisms,  
\begin{itemize}
\item[(i)] $\bar{1} \in (\YM(n,s)^{\Gamma})_{\h(n,s)}$ corresponds to $1 \in k \simeq H_{0}(\h(n,s))$,

\item[(ii)] $\overline{1 \otimes x_{2}}, \dots, \overline{1 \otimes x_{n}}, \overline{1 \otimes [x_{1},x_{2}]}, \dots, \overline{1 \otimes [x_{1},x_{n}]}, 
             \overline{1 \otimes z_{1}}, \dots, \overline{1 \otimes z_{s}} \in (\YM(n,s)^{\Gamma} \otimes \ym(n,s))_{\h(n,s)}$ 
            correspond to $q_{2}, \dots, q_{n}, p_{2}, \dots, p_{n}, z_{1}', \dots, z_{s}' \in U(n,s) \simeq H_{1}(\h(n,s),k)$,
\item[(ii)] $\overline{1 \otimes (\sum_{i=2}^{n} x_{i} \wedge [x_{1} , x_{i}] + (1/2) \sum_{a=1}^{s} z_{a} \wedge z_{a})} 
            \in (\YM(n,s)^{\Gamma} \otimes \Lambda^{2} \ym(n,s))_{\h(n,s)}$ 
            corresponds to the relation element given by $\sum_{i=2}^{n} [q_{i},p_{i}] + (1/2) \sum_{a=1}^{s} [z_{a}',z_{a}'] \in R(n,s)' \simeq H_{2}(\h(n,s))$.
\end{itemize}
The action of $\a = k.d$ on the homology groups is induced by the action of $x_{1}$ on the tensor component $\YM(n,s)^{\Gamma}$ 
of the modules $\YM(n,s)^{\Gamma} \otimes \Lambda^{\bullet} \ym(n,s)^{\Gamma}$ of the Chevalley-Eilenberg resolution. 
This immediately implies that the action of $x_{1}$ on $H_{0}(\h(n,s),k) \simeq k$ is trivial, for $x_{1} = d^{CE}_{1} (1 \otimes x_{1})$ is a boundary. 

Concerning the action of $\a$ on $H_{1}(\h(n,s),k)$, we have the identity $x_{1} . q_{i} = p_{i}$, for all $i = 2, \dots, n$, since 
\[     x_{1} \otimes x_{i} = 1 \otimes [x_{1},x_{i}] + x_{i} \otimes x_{1} + d^{CE}_{2}(1 \otimes (x_{1} \wedge x_{i})),      \]
and $\overline{x_{i} \otimes x_{1}} \in (\YM(n,s)^{\Gamma} \otimes \ym(n,s)^{\Gamma})_{\h(n,s)}$ vanishes. 
Moreover, the action of $x_{1}$ on both $p_{i}$ and $z_{a}$ is zero, for $i = 2, \dots, n$ and for $a = 1, \dots, s$, respectively. 
To prove the first half the statement, note that 
\begin{align*}
   x_{1} \otimes [x_{1},x_{i}] &= [x_{1},x_{i}] \otimes x_{1} + 1 \otimes [x_{1},[x_{1},x_{i}]] + d^{CE}_{2}(1 \otimes (x_{1} \wedge [x_{1},x_{i}]))
   \\                          
                               &= [x_{1},x_{i}] \otimes x_{1} - \sum_{j=2}^{n} 1 \otimes [x_{j},[x_{j},x_{i}]] 
                                + \frac{1}{2} \sum_{a,b=1}^{s} \Gamma_{a,b}^{i} 1 \otimes [z_{a},z_{b}] 
                                + d^{CE}_{2}(1 \otimes (x_{1} \wedge [x_{1},x_{i}]))
   \\ 
                               &= [x_{1},x_{i}] \otimes x_{1} + d^{CE}_{2}(\sum_{j=2}^{n} 1 \otimes x_{j} \wedge [x_{j},x_{i}]) 
                                - \sum_{j=2}^{n} (x_{j} \otimes [x_{j},x_{i}] - [x_{j},x_{i}] \otimes x_{j})
   \\
                               & - d^{CE}_{2}\Big(\frac{1}{2} \sum_{a,b=1}^{s} \Gamma_{a,b}^{i} 1 \otimes z_{a} \wedge z_{b}\Big)
                                 + \sum_{a,b=1}^{s} \Gamma_{a,b}^{i} z_{a} \otimes z_{b}
                                 + d^{CE}_{2}(1 \otimes (x_{1} \wedge [x_{1},x_{i}])),  
\end{align*}
and the elements $\overline{[x_{1},x_{i}] \otimes x_{1}}$, $\overline{x_{j} \otimes [x_{j},x_{i}] - [x_{j},x_{i}] \otimes x_{j}}$ 
and $\overline{z_{a} \otimes z_{b}}$ of the space $(\YM(n,s)^{\Gamma} \otimes \ym(n,s)^{\Gamma})_{\h(n,s)}$ vanish. 

The second half of the statement follows from 
\[     x_{1} \otimes z_{a} = 1 \otimes [x_{1},z_{a}] + z_{a} \otimes x_{1} + d^{CE}_{2}(1 \otimes (x_{1} \wedge z_{a}))     \]
and the fact that $\overline{z_{a} \otimes x_{1}} \in (\YM(n,s)^{\Gamma} \otimes \ym(n,s)^{\Gamma})_{\h(n,s)}$ vanishes. 

Finally, the action of $x_{1}$ on $H_{2}(\h(n,s),k) \simeq R(n,s)'$ is trivial, which can be proved as follows. 
First, note that 
\begin{align*}
   x_{1} \otimes \Big(\sum_{i=2}^{n} x_{i} \wedge [x_{1} , x_{i}] + \frac{1}{2} \sum_{a=1}^{s} z_{a} \wedge z_{a}\Big) 
                              &= d^{CE}_{3}(1 \otimes (x_{1} \wedge x_{i} \wedge [x_{1},x_{i}])) 
                               + \sum_{i=2}^{n} (x_{i} \otimes x_{1} \wedge [x_{1},x_{i}] - [x_{1},x_{i}] \otimes x_{1} \wedge x_{i})  
   \\
                              & - \sum_{i=2}^{n} 1 \otimes x_{1} \wedge [x_{i},[x_{1},x_{i}]] 
                                + \sum_{i=2}^{n} 1 \otimes x_{i} \wedge [x_{1},[x_{1},x_{i}]]
                                + \sum_{a=1}^{s} z_{a} \otimes x_{1} \wedge z_{a}
   \\                          
                              & + d^{CE}_{3}\Big(\frac{1}{2} \sum_{a=1}^{s} 1 \otimes x_{1} \wedge z_{a} \wedge z_{a}\Big) 
                                + \sum_{a=1}^{s} 1 \otimes [x_{1},z_{a}] \wedge z_{a}
                                - \frac{1}{2} \sum_{a=1}^{s} 1 \otimes x_{1} \wedge [z_{a}, z_{a}]
   \\
                              &= d^{CE}_{3}(1 \otimes (x_{1} \wedge x_{i} \wedge [x_{1},x_{i}])) 
                               + \sum_{i=2}^{n} (x_{i} \otimes x_{1} \wedge [x_{1},x_{i}] - [x_{1},x_{i}] \otimes x_{1} \wedge x_{i})
   \\
                              & + \sum_{i=2}^{n} 1 \otimes x_{i} \wedge [x_{1},[x_{1},x_{i}]]
                                + \sum_{a=1}^{s} z_{a} \otimes x_{1} \wedge z_{a}
                                + d^{CE}_{3}\Big(\frac{1}{2} \sum_{a=1}^{s} 1 \otimes x_{1} \wedge z_{a} \wedge z_{a}\Big)                           
   \\
                               & + \sum_{a=1}^{s} 1 \otimes [x_{1},z_{a}] \wedge z_{a}, 
\end{align*}
where we have simplified the last member using the first relation $r_{0,1}$ of \eqref{eq:relort}. 

Since the elements $\overline{x_{i} \otimes x_{1} \wedge [x_{1},x_{i}] - [x_{1},x_{i}] \otimes x_{1} \wedge x_{i}}$ and 
$\overline{z_{a} \otimes x_{1} \wedge z_{a}}$ of $(\YM(n,s)^{\Gamma} \otimes \Lambda^{2} \ym(n,s)^{\Gamma})_{\h(n,s)}$ vanish, 
it suffices to prove that the element of $(\YM(n,s)^{\Gamma} \otimes \Lambda^{2} \ym(n,s)^{\Gamma})_{\h(n,s)}$ induced by
\[    \sum_{i=2}^{n} 1 \otimes x_{i} \wedge [x_{1},[x_{1},x_{i}]] + \sum_{a=1}^{s} 1 \otimes [x_{1},z_{a}] \wedge z_{a}     \]
also vanish, which follows from 
\begin{align*}
   \sum_{i=2}^{n} 1 \otimes x_{i} \wedge [x_{1},[x_{1},x_{i}]] + \sum_{a=1}^{s} 1 \otimes [x_{1},z_{a}] \wedge z_{a} 
   &= - \sum_{i,j=2}^{n} 1 \otimes x_{i} \wedge [x_{j},[x_{j},x_{i}]] 
      + \frac{1}{2} \sum_{i=2}^{n} \sum_{a,b=1}^{s} \Gamma^{i}_{a,b} 1 \otimes x_{i} \wedge [z_{a},z_{b}] 
   \\   
   & - \sum_{a,b=1}^{s} \sum_{i=2}^{s} \Gamma^{i}_{a,b} 1 \otimes [x_{i},z_{b}] \wedge z_{a}
   \\
   &= \frac{1}{2} d^{CE}_{3}\Big(\sum_{i,j=2}^{n} 1 \otimes x_{i} \wedge x_{j} \wedge [x_{i},x_{j}] 
          + \sum_{i=2}^{n} \sum_{a,b=1}^{s} \Gamma^{i}_{a,b} 1 \otimes x_{i} \wedge z_{b} \wedge z_{a}\Big) 
   \\
   & - \frac{1}{2} \sum_{i=2}^{n} \sum_{a,b=1}^{s} \Gamma^{i}_{a,b} (x_{i} \otimes z_{b} \wedge z_{a} - 2 z_{b} \otimes x_{i} \wedge z_{a}) 
     - \sum_{i,j=2}^{n} x_{i} \otimes x_{j} \wedge [x_{i},x_{j}] 
   \\
   &- \frac{1}{2} \sum_{i,j=2}^{n} [x_{i},x_{j}] \otimes x_{i} \wedge x_{j}. 
\end{align*}

We now proceed to compute the homology groups $H_{\bullet}(\a,H_{\bullet}(\h(n,s),k))$. 
In this case, the Chevalley-Eilenberg resolution is of the form 
\[     0 \rightarrow \U(\a) \otimes k.d \rightarrow \U(\a) \overset{\epsilon_{\a}}{\rightarrow} k \rightarrow 0,     \]
where the first nonzero morphism is given by the restriction of the multiplication map of $\U(\a)$. 
In particular, $H_{0}(\a,k) \simeq k$, $H_{1}(\a,k) \simeq k[-2]$, $H_{0}(\a,R(n,s)') \simeq k[-6]$, $H_{1}(\a,R(n,s)') \simeq k[-8]$, and zero otherwise. 
Moreover, the previous resolution tells us that 
\begin{align*}
   H_{0}(\a,U(n,s)) &\simeq \mathrm{span}_{k} \cl{q_{2}, \dots, q_{n}, z_{1}', \dots, z_{s}'},
   \\ 
   H_{1}(\a,U(n,s)) &\simeq \mathrm{span}_{k} \cl{p_{2}, \dots, p_{n}, z_{1}', \dots, z_{s}'}[-2],
\end{align*} 
and the other homology groups vanish. 

Now, since the Hochschild-Serre spectral sequence $E_{p,q}^{2} = H_{p}(\a,H_{q}(\h(n,s),k)) \Rightarrow H_{p+q}(\ym(n,s)^{\Gamma},k)$ 
is concentrated in columns $p = 0 ,1$, it follows that $H_{n}(\ym(n,s)^{\Gamma},k) \simeq H_{0}(\a,H_{n}(\h(n,s),k)) \oplus H_{1}(\a,H_{n-1}(\h(n,s),k))$, 
for all $n \in \NN_{0}$ (see \cite{Rot09}, Cor. 10.29). 
Hence, 
\begin{equation}
\label{eq:homo}
     H_{\bullet}(\ym(n,s)^{\Gamma},k) \simeq \begin{cases}
                                                 k, &\text{if $\bullet = 0$,}
                                                 \\
                                                 V(n,s), &\text{if $\bullet = 1$,}
                                                 \\
                                                 R(n,s)^{\Gamma}, &\text{if $\bullet = 2$,} 
                                                 \\ 
                                                 k.\omega, &\text{if $\bullet = 3$,} 
                                                 \\
                                                 0,             &\text{else,}
                                           \end{cases}
\end{equation}
where $\omega = \sum_{i=1}^{n} x_{i} \otimes r_{0,i} + \sum_{a=1}^{s} z_{a} \otimes r_{1,a} \in (TV)_{8}$. 
Note that $k.\omega \simeq k[-8]$ and $R(n,s)^{\Gamma} \simeq V(n,s)^{*}[-8]$ in the category of graded $\so(n)$-modules provided with morphisms of degree zero. 
Also note that $\omega \in (V \otimes R) \cap (R \otimes V)$, because of the identity 
\[     \sum_{i=1}^{n} [x_{i}, r_{0,i}] + \sum_{a=1}^{s} [z_{a}, r_{1,a}] = 0     \]
in the tensor algebra $TV$. 

Moreover, we have the following: 
\begin{proposition}
\label{prop:projres}
Let $(n,s) \in (\NN \times \NN_{0}) \setminus \{ (1,0), (1,1) \}$.  
The minimal projective resolution $(K_{\bullet}'(\ym(n,s)^{\Gamma}),b_{\bullet}')$ of the $\YM(n,s)^{\Gamma}$-module $k$ is given by 
\begin{equation}
\label{eq:comp}
     0 \rightarrow \YM(n,s)^{\Gamma}[-8] \overset{b_{3}'}{\rightarrow} \YM(n,s)^{\Gamma} \otimes R(n,s)^{\Gamma} \overset{b_{2}'}{\rightarrow}
                     \YM(n,s)^{\Gamma} \otimes V(n,s) \overset{b_{1}'}{\rightarrow} \YM(n,s)^{\Gamma} \overset{b_{0}'}{\rightarrow} k \rightarrow 0,     
\end{equation}
with differential  
\begin{equation}
\label{eq:dif}
\begin{split}
   b_{3}'(z) &= \sum_{i=1}^{n} z x_{i} \otimes r_{0,i} + \sum_{a=1}^{s} z z_{a} \otimes r_{1,a}, 
   \\
   b_{2}'(z \otimes r_{0,i}) &= \sum_{j=1}^{n} (z x_{j}^{2} \otimes x_{i} - 2 z x_{j} x_{i} \otimes x_{j} + z x_{i} x_{j} \otimes x_{j}) 
                               - \sum_{a,b=1}^{s} \Gamma_{a,b}^{i} z z_{a} \otimes z_{b},
   \\
   b_{2}'(z \otimes r_{1,a}) &= \sum_{i=1}^{n} \sum_{b=1}^{s} \Gamma^{i}_{a,b} (z x_{i} \otimes z_{b} - z z_{b} \otimes x_{i}),
   \\
   b_{1}'(z \otimes x_{i}) &= z x_{i},
   \\
   b_{1}'(z \otimes z_{a}) &= z z_{a}, 
   \\
   b_{0}'(z) &= \epsilon_{\ym(n,s)^{\Gamma}}(z).
\end{split}
\end{equation}
\end{proposition}
\noindent\textbf{Proof.}
As stated above, the $\YM(n,s)^{\Gamma}$-modules $K_{\bullet}'(\ym(n,s)^{\Gamma},k)$ giving the minimal projective resolution of $k$ are 
of the form $\YM(n,s)^{\Gamma} \otimes_{k} \Tor^{\YM(n,s)^{\Gamma}}_{\bullet}(k,k)$. 
Since $H_{\bullet}(\ym(n,s)^{\Gamma},k) \simeq \Tor^{\YM(n,s)^{\Gamma}}_{\bullet}(k,k)$, the previous homological computations tell us that the modules involved 
in \eqref{eq:comp} are correct and we only need to prove that the claimed differential provides a resolution of $k$. 
It is obvious to see that this gives a complex, \textit{i.e.} $b_{i}' \circ b_{i+1}' = 0$, for $i = 0, 1, 2$. 

A similar argument to the one given after Lemma 2.3 in \cite{BG06} shows that  
\[     \YM(n,s)^{\Gamma} \otimes R(n,s)^{\Gamma} \overset{b_{2}'}{\rightarrow}
                     \YM(n,s)^{\Gamma} \otimes V(n,s) \overset{b_{1}'}{\rightarrow} \YM(n,s)^{\Gamma} \overset{b_{0}'}{\rightarrow} k \rightarrow 0,     \]
is exact. 
Indeed, let $A = TV/\cl{R}$ be a connected graded algebra, where $V = \oplus_{i \in \NN} V_{i}$ and $R = \oplus_{i \in \NN_{\geq 2}} R_{i} \subseteq T^{+}V = \oplus_{i \in \NN} V^{\otimes i}$ are 
nonnegatively graded vector spaces and we further suppose $R \cap ((T^{+}V) R (TV) + (TV) R (T^{+}V)) = 0$, in order to avoid ambiguities. 
In this case, we consider the complex 
\[     A \otimes R \overset{b_{2}'}{\rightarrow} A \otimes V \overset{b_{1}'}{\rightarrow} A \overset{b_{0}'}{\rightarrow} k \rightarrow 0,     \]
where $b_{0}'$ is the augmentation of $A$, $b_{1}'$ is the restriction of the multiplication of $A$, and $b_{2}'$ is the $A$-linear extension 
of the morphism $R \rightarrow TV \otimes V \rightarrow A \otimes V$, given by the composition of the inclusion of $R$ in $T^{+}V \simeq TV \otimes V$ with the canonical surjection. 
Since the zeroth degree component of $b_{1}'$ is zero, whereas the $m$-th degree component ($m \geq 1$) is the surjective map 
\[     \sum_{i \in \NN} \frac{(TV)_{m-i} \otimes V_{i}}{\cl{R}_{m-i} V_{i}} \rightarrow \frac{(TV)_{m}}{\cl{R}_{m}},     \]
where $\cl{R}_{m}$ is the component of degree $m$ of the ideal generated by $R$ in $TV$, 
it follows that $\Ker(b_{0}') = \mathrm{Im}(b_{1}')$. 
Moreover, taking into account that the $m$-th degree component of $b_{2}'$ is the map
\[     \frac{\sum_{i \in \NN} (TV)_{m-i} R_{i}}{\sum_{i \in \NN} \cl{R}_{m-i} R_{i}} 
       \rightarrow \frac{\sum_{j \in \NN} (TV)_{m-j} V_{j}}{\sum_{j \in \NN} \cl{R}_{m-j} V_{j}},     \]
then 
\[     \mathrm{Im}(b_{2}')_{m} 
             = \frac{\sum_{i \in \NN} (TV)_{m-i} R_{i} + \sum_{j \in \NN} \cl{R}_{m-j} V_{j}}{\sum_{j \in \NN} \cl{R}_{m-j} V_{j}},     \]
which coincides with 
\[     \Ker(b_{1}')_{m} = \frac{\cl{R}_{m}}{\sum_{j \in \NN} \cl{R}_{m-j} V_{j}}.     \]

It is still left to show that $b_{3}'$ is injective and that $\mathrm{Im}(b_{3}') = \Ker(b_{2}')$. 
Let $z \in \YM(n,s)^{\Gamma}$ such that $b_{3}'(z) = 0$, \textit{i.e.} such that $z x_{i} = 0$ and $z z_{a} = 0$ for all $i = 1, \dots, n$ and $a = 1, \dots, s$. 
Since $\YM(n,s)^{\Gamma}$ is a free right module over the integral domain $\U(\ym(n,s)^{\Gamma}_{0})$, we conclude that $z = 0$, so $b_{3}'$ is injective. 

In oder to prove that $\mathrm{Im}(b_{3}') = \Ker(b_{2}')$, we proceed as follows. 
Since the third component of the minimal projective resolution of $k$ is $\YM(n,s)^{\Gamma} \otimes k[-8]$, 
the kernel of $b_{2}'$ should be of the form $i : \YM(n,s)^{\Gamma} \otimes k[-8] \rightarrow \YM(n,s)^{\Gamma} \otimes R(n,s)^{\Gamma}$. 
Hence, we have the commutative diagram 
\[
\xymatrix@C-10pt
{
\YM(n,s)^{\Gamma} \otimes k[-8]
\ar[rd]^{b_{3}'}
\ar@{-->}[dd]^{\bar{i}}
&
&
\\
&
\YM(n,s)^{\Gamma} \otimes R(n,s)^{\Gamma}
\ar[r]^{b_{2}'}
&
\YM(n,s)^{\Gamma} \otimes V(n,s)
\\
\YM(n,s)^{\Gamma} \otimes k[-8]
\ar[ru]^{i}
&
&
}
\] 
Since $b_{3}'$ is injective, $\bar{i}$ is also so, which in turn implies that it is an isomorphism, for it is an injective endomorphism (of degree zero) 
of a locally finite dimensional graded vector space. 
Hence, $\mathrm{Im}(b_{3}') = \Ker(b_{2}')$, as was to be proved. 
\qed

\begin{remark}
Note that for $(n,s) = (1,1)$ the previous complex does not yield a resolution of $k$, since $H_{\bullet}(\h(1,1),k) \simeq k$, for $\bullet \in \NN_{0}$.                        
\end{remark}

There is an analogous free resolution of right $\YM(n,s)^{\Gamma}$-modules of the trivial module $k$, for $(n,s) \neq (1,0), (1,1)$, which is of the form 
\begin{equation}
\label{eq:comp*}
     0 \rightarrow \YM(n,s)^{\Gamma}[-8] \overset{b_{3}''}{\rightarrow} R(n,s)^{\Gamma} \otimes \YM(n,s)^{\Gamma} \overset{b_{2}''}{\rightarrow}
                     V(n,s) \otimes \YM(n,s)^{\Gamma} \overset{b_{1}''}{\rightarrow} \YM(n,s)^{\Gamma} \overset{b_{0}'}{\rightarrow} k \rightarrow 0,     
\end{equation}
with differential given by $b_{0}'' = \epsilon_{\ym(n,s)^{\Gamma}}$ and  
\begin{equation}
\label{eq:dif*}
\begin{split}
   b_{3}''(z) &= \sum_{i=1}^{n} r_{0,i} \otimes x_{i} z  - \sum_{a=1}^{s} r_{1,a} \otimes z_{a} z, 
   \\
   b_{2}''(r_{0,a} \otimes z) &= \sum_{j=1}^{n} (x_{i} \otimes x_{j}^{2} z - 2 x_{j} \otimes x_{i} x_{j} z + x_{j} \otimes x_{j} x_{i} z) 
                               - \sum_{a,b=1}^{s} \Gamma_{a,b}^{i} z_{a} \otimes z_{b} z,
   \\
   b_{2}''(r_{1,a} \otimes z) &= \sum_{i=1}^{n} \sum_{b=1}^{s} \Gamma^{i}_{a,b} (x_{i} \otimes z_{b} z - z_{b} \otimes x_{i} z),
   \\
   b_{1}''(x_{i} \otimes z) &= x_{i} z,
   \\
   b_{1}''(z_{a} \otimes z) &= z_{a} z.  
\end{split}
\end{equation}
We shall denote this resolution by $(K_{\bullet}''(\ym(n,s)^{\Gamma}),b_{\bullet}'')$. 

We now recall that a nonnegatively graded connected algebra $A$ is called \emph{(left) AS-regular of dimension $d$ and of Gorenstein parameter $l$} 
if it has finite (left) global dimension $d$ and it satisfies that 
\[     \mathcal{E}xt^{i}_{A} (k, A) \simeq \begin{cases}
                                                                              k[l], &\text{if $i = d$},
                                                                              \\
                                                                              0, &\text{else}. 
                                                                           \end{cases}
\]
We point out that no noetherianity assumption on $A$ is required (see \cite{MM11}, Def. 1.1). 
Since we shall be dealing with graded Hopf algebras, both left and right definitions of AS-regular coincide in this case, 
so we will omit the reference to the side. 

Using the fact that $\mathcal{H}om_{\YM(n,s)^{\Gamma}}(\YM(n,s)^{\Gamma}[i],\place) \simeq (\place)[-i]$, that 
\begin{align*}
     \YM(n,s)^{\Gamma} \otimes V(n,s) &\simeq (\YM(n,s)^{\Gamma})^{n}[-2] \oplus (\YM(n,s)^{\Gamma})^{s}[-3],
     \\ 
     \YM(n,s)^{\Gamma} \otimes R(n,s)^{\Gamma} &\simeq (\YM(n,s)^{\Gamma})^{n}[-6] \oplus (\YM(n,s)^{\Gamma})^{s}[-5],
\end{align*} 
and an elementary computation involving the differentials, 
one sees that the functor $\mathcal{H}om_{\YM(n,s)^{\Gamma}}(\place,\YM(n,s)^{\Gamma})$ sends the resolution 
$(K'_{\bullet}(\ym(n,s)^{\Gamma}),b_{\bullet}')$ 
of left $\YM(n,s)^{\Gamma}$-modules of $k$ to the shift $(K''_{\bullet}(\ym(n,s)^{\Gamma}),b_{\bullet}'')[8]$. 
We remark that the we are always using Koszul's sign rule, and, in particular, 
$X \simeq X^{**}$, via the map $x \mapsto (f \mapsto (-1)^{|x||f|}) f(x)$; 
$X \otimes Y \simeq Y \otimes X$, via $x \otimes y \mapsto (-1)^{|y||x|} y \otimes x$, and 
$Y \otimes X^{*} \simeq \mathcal{H}om_{k}(X,Y)$, via $y \otimes f \mapsto (x \mapsto y f(x))$, 
for (finite dimensional) graded (resp. super) vector spaces $X$ and $Y$, and  homogeneous elements $x \in X$, $y \in Y$ and $f \in X^{*}$. 
This tells us that $\YM(n,s)^{\Gamma}$ is AS-regular with Gorenstein parameter $8$, in the sense of \cite{MM11}, for $(n,s) \neq (1,1)$. 

As a corollary of the description of the left and right $\YM(n,s)^{\Gamma}$-module resolutions of $k$, 
we may obtain the minimal projective resolution of the $\YM(n,s)^{\Gamma}$-bimodule $\YM(n,s)^{\Gamma}$, for $(n,s) \neq (1,0), (1,1)$. 
Consider first the free graded $\YM(n,s)^{\Gamma}$-bimodule 
$K_{\bullet}(\ym(n,s)^{\Gamma}) = \YM(n,s)^{\Gamma} \otimes H_{\bullet}(\ym(n,s)^{\Gamma},k) \otimes \YM(n,s)^{\Gamma}$, for $\bullet \geq 0$, 
given with the obvious action, where $H_{\bullet}(\ym(n,s)^{\Gamma},k)$ is identified to a subspace of $TV(n,s)$ following \eqref{eq:homo}. 
Then, we define on it the differential given by $b_{1} = b_{1}' \otimes 1_{\YM(n,s)^{\Gamma}} - 1_{\YM(n,s)^{\Gamma}} \otimes b_{1}''$, 
$b_{3} = b_{3}' \otimes 1_{\YM(n,s)^{\Gamma}} - 1_{\YM(n,s)^{\Gamma}} \otimes b_{3}''$, and 
\[     b_{2} (y \otimes y_{1} \dots y_{r} \otimes y') = \sum_{i=1}^{r} y y_{1} \dots y_{i-1} \otimes y_{i} \otimes y_{i+1} \dots y_{r} y',     \]
where $y, y' \in \YM(n,s)^{\Gamma}$, $y_{1} \dots y_{r} \in R(n,s)^{\Gamma}$, for $y_{i} \in V(n,s)$, $i = 1, \dots, r$. 
Set $b_{\bullet} = 0$, for $\bullet \geq 4$. 
It is easily proved to be a complex. 
If we define $b_{0} : K_{\bullet}(\ym(n,s)^{\Gamma}) \rightarrow \YM(n,s)^{\Gamma}$ given by the multiplication, we obtain an augmented complex.   
Moreover, taking into account the isomorphism 
$(K_{\bullet}(\ym(n,s)^{\Gamma}), b_{\bullet}) \otimes_{\YM(n,s)^{\Gamma}} k \simeq (K'_{\bullet}(\ym(n,s)^{\Gamma}),b_{\bullet}')$ 
of complexes of left $\YM(n,s)^{\Gamma}$-modules, 
\cite{BM06}, Prop. 4.1, tells us that the bimodule complex $(K_{\bullet}(\ym(n,s)^{\Gamma}), b_{\bullet})$ is in fact a resolution of $\YM(n,s)^{\Gamma}$. 

Let us recall that a locally finite dimensional nonnegatively graded and connected algebra $A$ is called 
\emph{(left) graded Calabi-Yau of dimension $d$ and of parameter $l$} 
if it has a finite resolution composed of finitely generated projective bimodules, finite global dimension $d$, and 
satisfies that (\textit{cf.} \cite{BT07}, Def. 4.2)
\[     \mathcal{E}xt^{i}_{A^{e}} (A, A^{e}) \simeq \begin{cases}
                                                                                                           A[l], &\text{if $i = d$},
                                                                                                           \\
                                                                                                           0, &\text{else}, 
                                                                                                       \end{cases}
\]
in the category of (right) $A^{e}$-modules. 
Again, taking into account that we shall be dealing with graded Hopf algebras, both left and right definitions of Calabi-Yau coincide, 
and we will thus omit the reference to the side. 
 
We claim that the image of $(K_{\bullet}(\ym(n,s)^{\Gamma}), b_{\bullet})$ 
under the functor $\mathcal{H}om_{(\YM(n,s)^{\Gamma})^{e}}(\place, (\YM(n,s)^{\Gamma})^{e})$ 
is isomorphic to its shift $(K_{\bullet}(\ym(n,s)^{\Gamma}), (-1)^{\bullet} b_{\bullet})[8]$, 
in the category of graded $\YM(n,s)^{\Gamma}$-bimodules, which 
follows easily as before from the fact that 
\[     \mathcal{H}om_{(\YM(n,s)^{\Gamma})^{e}}(\YM(n,s)^{\Gamma} \otimes X \otimes \YM(n,s)^{\Gamma}, (\YM(n,s)^{\Gamma})^{e}) \simeq 
       \mathcal{H}om_{k}(X, (\YM(n,s)^{\Gamma})^{e}) \simeq \YM(n,s)^{\Gamma} \otimes X^{*} \otimes \YM(n,s)^{\Gamma},     \]
where $X$ is a finite dimensional graded vector space, and an elementary computation with the differential. 
We remark that, if $X = \oplus_{i \in \ZZ} X_{i}$, then $X^{*} = \oplus_{i \in \ZZ} (X^{*})_{i}$, where $(X^{*})_{i} = (X_{-i})^{*}$. 
This in turn implies that the super Yang-Mills algebra is graded Calabi-Yau of dimension $3$ and of parameter $8$, for $(n,s) \neq (1,0), (1,1)$. 

Hence, we have obtained the following: 
\begin{proposition}
Let $(n,s) \in (\NN \times \NN_{0}) \setminus \{ (1,0), (1,1) \}$. 
The graded algebra $\YM(n,s)^{\Gamma}$ is AS-regular of global dimension $3$ and of Gorenstein parameter $8$. 
Moreover, it is also graded Calabi-Yau of the same dimension and of the same parameter. 
\end{proposition}

\subsection{Some consequences}
\label{subsec:cons}

\subsubsection{The Hilbert series of the super Yang-Mills algebra}

We shall compute the Hilbert series of both the graded algebra $\YM(n,s)^{\Gamma}$ and the graded Lie algebra $\ym(n,s)^{\Gamma}$. 
In order to do so, we first recall that the Hilbert series is an Euler-Poincar\'e map in the category of graded vector spaces (see \cite{Lang02}, Chap. III, \S 8) 
and that the Euler-Poincar\'e characteristic of a complex of graded vector spaces, with respect to the Euler-Poincar\'e map given by taking Hilbert series, 
coincides with that of its homology (see \cite{Lang02}, Chap. XX, \S 3, Thm. 3.1). 

\begin{proposition}
The Hilbert series of the graded algebra $\YM(n,s)^{\Gamma}$ provided with the grading explained in Subsection \ref{subsec:defi}, \textit{i.e.} 
such that $x_{i}$ has degree $2$, for all $i = 1, \dots, n$, and $z_{a}$ has degree $3$, for all $a = 1, \dots, s$, is given by
\[     \YM(n,s)^{\Gamma}(t) = \frac{1}{1 - n t^{2} - s t^{3} + s t^{5} + n t^{6} - t^{8}}.     \]
\end{proposition}
\noindent\textbf{Proof.}
Since the complex \eqref{eq:comp} of graded vector spaces is exact, its Euler-Poincar\'e characteristic, 
with respect to the Euler-Poincar\'e map given by taking Hilbert series, coincides with that of its homology. 
Hence, 
\[     1 - \YM(n,s)^{\Gamma}(t) + (n t^{2} + s t^{3}) \YM(n,s)^{\Gamma}(t) - (n t^{6} + s t^{5}) \YM(n,s)^{\Gamma}(t) + t^{8} \YM(n,s)^{\Gamma}(t) = 0,     \]
from which the proposition follows. 
\qed

We state the following result, which we believe is well-known, but we give its proof because we do not know any specific reference. 
\begin{proposition}
Let $V$ be a positively graded vector space with Hilbert series given by $\sum_{j \in \NN} \nu_{j} t^{j}$. 
Then, 
\[     \nu_{j} = \frac{(-1)^{j}}{j} \sum_{d|j} (-1)^{d} a_{d} \mu(\frac{j}{d}),     \]
for all $j \in \NN$, where $\mu$ is the M\"obius function and the coefficients $a_{d}$ are obtained from the formal series 
\[     \log(S(V)(t)) = \sum_{d \in \NN} \frac{a_{d}}{d} t^{d}.     \]
Equivalently, if $S(V)(t)^{-1}$ is a polynomial of degree $m$ with roots $\lambda_{i}$, $i = 1, \dots, m$, 
\[     a_{d} = \sum_{i=1}^{m} \lambda_{i}^{-d}.     \]
\end{proposition}
\noindent\textbf{Proof.}
First note that the Hilbert series of the symmetric algebra (in the super sense) of a graded vector space $V$ satisfying $V(t) = \sum_{j \in \NN} \nu_{j} t^{j}$ 
is given by 
\[     S(V)(t) = \frac{\prod_{i \in 2 \NN - 1} (1 + t^{i})^{\nu_{i}}}{\prod_{i \in 2 \NN} (1 - t^{i})^{\nu_{i}}} 
               = \prod_{i \in \NN} (1 - (-1)^{i} t^{i})^{-(-1)^{i} \nu_{i}}.     \]
In consequence, using that $\log(1-t) = - \sum_{i \in \NN} t^{i}/i$, we obtain that 
\[     \log(S(V)(t)) = - \sum_{i \in \NN} (-1)^{i} \nu_{i} \log(1 - (-1)^{i} t^{i}) = \sum_{i,j \in \NN} (-1)^{i} \nu_{i} \frac{(-t)^{i.j}}{j} 
                     = \sum_{l \in \NN} \sum_{i|l} (-1)^{i+l} \nu_{i} i \frac{t^{l}}{l},     \]
which means that 
\[     a_{d} = (-1)^{d} \sum_{i|d} (-1)^{i} \nu_{i} i.     \]
By the M\"obius inversion formula we have thus 
\[     \nu_{j} = \frac{(-1)^{j}}{j} \sum_{d|j} (-1)^{d} a_{d} \mu(\frac{j}{d}),     \]
as was to be shown. 

To prove the last statement, note that, if $S(V)(t)^{-1}$ is a polynomial of the form 
\[     S(V)(t)^{-1} = \prod_{i=1}^{m} (1 - \frac{t}{\lambda_{i}}),    \]
then 
\[     \log(S(V)(t)) = \sum_{i=1}^{m} \sum_{j \in \NN} \frac{1}{j} \Big(\frac{t}{\lambda_{i}}\Big)^{j},     \]
and the proposition is proved. 
\qed

As a consequence of the two previous propositions and the PBW theorem for graded Lie algebras we obtain the
\begin{corollary}
\label{cor:hilserym}
Let $\sum_{j \in \NN} \nu(n,s)_{j} t^{j}$ be the Hilbert series of the graded Lie algebra $\ym(n,s)^{\Gamma}$ provided with the grading explained in Subsection \ref{subsec:defi}. 
We shall sometimes write $\nu_{j}$ instead of $\nu(n,s)_{j}$. 
Then, 
\[     \nu_{j} = \frac{(-1)^{j}}{j} \sum_{d|j} (-1)^{d} a_{d} \mu(\frac{j}{d}),     \]
for all $j \in \NN$, where the coefficients $a_{d}$ are obtained from the formal series 
\[     \log(1 - n t^{2} - s t^{3} + s t^{5} + n t^{6} - t^{8}) = - \sum_{d \in \NN} \frac{a_{d}}{d} t^{d},     \]
or equivalently, if $\lambda_{1}, \dots, \lambda_{8}$ are the roots of $1 - n t^{2} - s t^{3} + s t^{5} + n t^{6} - t^{8}$, we have that 
\[     a_{d} = \sum_{i=1}^{8} \lambda_{i}^{-d} = (-1)^{d} \sum_{i=1}^{8} \lambda_{i}^{d}.     \]
\end{corollary}

\subsubsection{Some explicit computations of basis elements of a super Yang-Mills algebra}

As a corollary of the previous computation of the Hilbert series of the super Yang-Mills algebra, we will present in this paragraph 
some simple calculations of the basis elements of $\ym(3,1)^{\Gamma}$ for low degree homogeneous components. 
Note that $\ym(3,1)^{\Gamma}$ cannot be equivariant. 

We remark that in this case, under the assumptions explained in Subsection \ref{subsec:defi}, the space of relations can be taken to be of the form 
\begin{align*}
  r_{0,1} &= [x_{2},[x_{2},x_{1}]] + [x_{3},[x_{3},x_{1}]] - \frac{1}{2} [z_{1},z_{1}],     
  \\
  r_{0,2} &= [x_{1},[x_{1},x_{2}]] + [x_{3},[x_{3},x_{2}]],
  \\
  r_{0,3} &= [x_{1},[x_{1},x_{3}]] + [x_{2},[x_{2},x_{3}]],
  \\
  r_{1,1} &= [x_{1},z_{1}].  
\end{align*}
We shall first consider the super Yang-Mills algebra $\ym(3,1)^{\Gamma}$ as a graded Lie algebra, and consider the descending sequence of graded ideals $\{ F^{\bullet}\ym(3,1)^{\Gamma} \}_{\bullet \in \NN_{0}}$, 
where $F^{j}\ym(3,1)^{\Gamma}$ is the graded vector subspace of $\ym(3,1)^{\Gamma}$ given by the elements of degree greater than or equal to $j+2$. 
Our aim is to give explicit bases for the quotients $\ym(3,1)^{\Gamma}/F^{j}\ym(3,1)^{\Gamma}$, for $j = 1, \dots, 7$.

Using Corollary \ref{cor:hilserym} for $\ym(3,1)^{\Gamma}$, the sequence of dimensions $\nu(3,1)_{j}$ ($j \in \NN$) is given by 
\[     0, 3, 1, 3, 2, 6, 6, 12, 15, 33, 42, 77, 114, 213, 314, 555, 876, 1540, 2460, 4242,\dots     \]

An ordered basis for the quotient algebra $\ym(3,1)^{\Gamma}/F^{1}\ym(3,1)^{\Gamma}$ of the super Yang-Mills algebra, which is concentrated in degree $2$, 
is given by $\B_{1} = \{ x_{1}, x_{2}, x_{3} \}$, 
whereas a basis of $\ym(3,1)^{\Gamma}/F^{2}\ym(3,1)^{\Gamma}$ may be defined as $\B_{2} = \{ x_{1}, x_{2}, x_{3}, z_{1} \}$. 
Notice that $\ym(3,1)^{\Gamma}/\C^{2}(\ym(3,1)^{\Gamma}) = \ym(3,1)^{\Gamma}/F^{2}\ym(3,1)^{\Gamma}$. 

For the quotient $\ym(3,1)^{\Gamma}/F^{3}\ym(3,1)^{\Gamma}$, a possible ordered basis is
\[     \B_{3} = \{ x_{1}, x_{2}, x_{3}, z_{1}, x_{12}, x_{13} , x_{23} \},     \]
where $x_{ij} = [x_{i} , x_{j}]$, ($i, j = 1, 2, 3$). 
To prove that it is indeed a basis we must only show that it generates $\ym(3,1)^{\Gamma}/F^{3}\ym(3,1)^{\Gamma}$ (since $\#(\B_{3}) = 7$),
which is obtained using that $x_{ij} = - x_{ji}$. 

The quotient $\ym(3,1)^{\Gamma}/F^{4}\ym(3,1)^{\Gamma}$ is a little more interesting, and a possible ordered basis for it is of the form 
\[     \B_{4} = \{ x_{1}, x_{2}, x_{3}, z_{1}, x_{12}, x_{13} , x_{23}, y_{2}, y_{3} \},     \]
where we denote $y_{i} = [x_{i},z_{1}]$ ($i = 1, 2, 3$). 
Again, it suffices to show that it generates $\ym(3,1)^{\Gamma}/F^{4}\ym(3,1)^{\Gamma}$, for $\#(\B_{4}) = 9$. 
But this follows from the fact that $y_{1} = [x_{1},z_{1}] = r_{1,1} = 0$.  
Note that $\ym(3,1)^{\Gamma}/\C^{3}(\ym(3,1)^{\Gamma}) = \ym(3,1)^{\Gamma}/F^{4}\ym(3,1)^{\Gamma}$, because the vanishing of $r_{0,1}$ implies that $[z_{1},z_{1}] = 0$. 

We claim that the following set is a basis of $\ym(3,1)^{\Gamma}/F^{5}\ym(3,1)^{\Gamma}$:
\[     \B_{5} = \{ x_{1}, x_{2}, x_{3}, z_{1}, x_{12}, x_{13}, x_{23}, y_{2}, y_{3}, x_{112}, x_{221}, x_{113}, x_{123}, x_{312}, [z_{1},z_{1}]  \},     \]
where $x_{ijk} = [x_{i},[x_{j},x_{k}]]$.
Indeed, as before, we only have to prove that the previous set it is a system of generators of $\ym(3,1)^{\Gamma}/F^{5}\ym(3,1)^{\Gamma}$, because $\#(\B_{5}) = 15$.
This is direct, and can be seen from the Yang-Mills relations
\[     x_{332} = - x_{112}, \hskip 0.5cm x_{331} = - x_{221} + \frac{1}{2} [z_{1},z_{1}], \hskip 0.5cm x_{223} = - x_{113},    \]
and the relations given by antisymmetry and Jacobi identity, \textit{i.e.} $x_{ijk} = - x_{ikj}$, and $x_{213} = x_{123} + x_{312}$. 

Now, a basis of $\ym(3,1)^{\Gamma}/F^{6}\ym(3,1)^{\Gamma}$ can be given by
\[     \B_{6} = \{ x_{1}, x_{2}, x_{3}, z_{1}, x_{12}, x_{13}, x_{23}, y_{2}, y_{3}, x_{112}, x_{221}, x_{113}, x_{123}, x_{312}, [z_{1},z_{1}], y_{12}, y_{13}, y_{22}, y_{23}, y_{32}, y_{33}  \},     \]
where we write $y_{ij} = [x_{i},[x_{j},z_{1}]]$. 
That this is a basis is immediate, since its cardinality is $21$ and it is a set of generators of $\ym(3,1)^{\Gamma}/F^{6}\ym(3,1)^{\Gamma}$, for $[z_{1},x_{ij}] = y_{ji} - y_{ij}$. 

The case of $\ym(3,1)^{\Gamma}/F^{7}\ym(3,1)^{\Gamma}$ is a little more complicated. 
We shall prove that
\begin{align*}
   \B_{7} = \{ &x_{1}, x_{2}, x_{3}, z_{1}, x_{12}, x_{13} , x_{23}, y_{2}, y_{3}, x_{112}, x_{221}, x_{113}, x_{123}, x_{312}, [z_{1},z_{1}], y_{12}, y_{13}, y_{22}, y_{23}, y_{32}, y_{33}
   \\
   &x_{1112}, x_{1221}, x_{1113}, x_{1123}, x_{2221}, x_{2113}, x_{2312}, x_{3112}, x_{3221}, x_{3312}, [x_{2},[z_{1},z_{1}]], [x_{3},[z_{1},z_{1}]] \},
\end{align*}
is an ordered basis, where $x_{ijkl} = [x_{i},[x_{j},[x_{k},x_{l}]]]$.

In order to prove that $\B_{7}$ is a basis it suffices again to verify that it is a system of generators of the super vector space 
$\ym(3,1)^{\Gamma}/F^{7}\ym(3,1)^{\Gamma}$.
On the one hand, taking into account that
\[     [[x_{i},x_{j}],[x_{k},x_{l}]] = [[[x_{i},x_{j}],x_{k}],x_{l}] + [x_{k},[[x_{i},x_{j}],x_{l}]]
                                     = [x_{l},[x_{k},[x_{i},x_{j}]]] + [x_{k},[x_{l},[x_{j},x_{i}]]]
                                     = x_{lkij} + x_{klji},
\]
that $[z_{1},[x_{i},z_{1}]] = [x_{i},[z_{1},z_{1}]]/2$, and that $[[[x_{i},x_{j}],x_{k}],x_{l}] = [x_{l},[x_{k},[x_{i},x_{j}]]] = x_{lkij}$, 
we see that the set $\B_{6} \cup \{ x_{ijkl} : i, j, k, l = 1, 2, 3 \} \cup \{ [x_{2},[z_{1},z_{1}]], [x_{3},[z_{1},z_{1}]] \}$ is a system of generators.
We shall prove that it is generated by $\B_{7}$.
In fact, we only need to prove that the latter generates the set
\[     \{  x_{i112}, x_{i221}, x_{i113}, x_{i123}, x_{i312} : i = 1, 2, 3 \},     \]
because $\{ x_{112}, x_{221}, x_{113}, x_{123}, x_{312} \}$ are generators of the homogeneous elements of degree $6$.
This last statement is direct:
\begin{align*}
     x_{3113} &= x_{1221}, \hskip 0.5cm x_{2112} = - x_{1221}, \hskip 0.5cm
     x_{2123} = x_{3221} + x_{2312} - x_{1113},
     \\
     x_{1312} &= \frac{x_{3112} + x_{2113} - x_{1123}}{2}, \hskip 0.95cm x_{3123} = \frac{x_{1112}+x_{2221}-x_{3312}}{2}- \frac{[x_{2},[z_{1},z_{1}]]}{4}.
\end{align*}

\subsubsection{\texorpdfstring{The ideal $\hat{\tym}(n,s)^{\Gamma}$ and some algebraic properties of the super Yang-Mills algebra}{The ideal tym and some algebraic properties of the super Yang-Mills algebra}}
\label{subsec:tym}

In this paragraph we shall obtain an important result which will be useful in the sequel: the graded Lie algebra $\hat{\tym}(n,s)^{\Gamma}$ is free. 
We will also compute the Hilbert series of its space of generators. 
In order to do so, we shall prove that the homology groups of degree greater than one of the graded Lie algebra $\hat{\tym}(n,s)^{\Gamma}$ with coefficients in $k$ vanish. 
At the end, we shall derive several algebraic properties of the super Yang-Mills algebras. 

First, we need a subsidiary result, for which we recall that the quotient $\ym(n,s)^{\Gamma}/\hat{\tym}(n,s)^{\Gamma} \simeq V(2)_{0}$ tells us 
that $V(2)_{0}$, and hence $S(V(2)_{0})$, is a module over $\YM(n,s)^{\Gamma}$. 
\begin{proposition}
\label{prop:homow}
Let $n \geq 2$. 
The homology groups $H_{\bullet} (\ym(n,s)^{\Gamma}, S(V(2)_{0}))$ of the graded Lie algebra $\ym(n,s)^{\Gamma}$ with coefficients in 
the module $\U(\ym(n,s)^{\Gamma}/\hat{\tym}(n,s)^{\Gamma}) \simeq S(V(2)_{0})$, which is obviously a graded right $\YM(n,s)^{\Gamma}$-module, are given by 
\[     H_{\bullet}(\ym(n,s)^{\Gamma},S(V(2)_{0})) \simeq \begin{cases}
                                                   k, &\text{if $\bullet = 0$,}
                                                   \\
                                                   \hat{W}(n,s)^{\Gamma}, &\text{if $\bullet = 1$,}
                                                   \\
                                                   0,             &\text{else,}
                                                   \end{cases}
\]
where $\hat{W}(n,s)^{\Gamma}$ is the graded vector space with Hilbert series 
\[     \hat{W}(n,s)^{\Gamma}(t) = (n-2) t^{2} + (2n-3) t^{4} + \sum_{k \geq 3} (2n - 4) t^{2 k} + \sum_{k \geq 1} s t^{2 k + 1}.     \]
\end{proposition}
\noindent\textbf{Proof.}
We shall make the computation of the homology groups using the minimal projective resolution in Proposition \ref{prop:projres}.  
It is direct to see that the complex $(S(V(2)_{0}) \otimes_{\YM(n,s)^{\Gamma}} K_{\bullet}(\ym(n,s)^{\Gamma}), 1 \otimes b_{\bullet}')$ is the direct sum 
of two complexes $(K_{\bullet}^{\sim},b_{\bullet}^{\sim})$ and $(K_{\bullet}^{\backsim},b_{\bullet}^{\backsim})$, which we now describe. 
The latter is given by $K_{2}^{\backsim} = S(V(2)_{0}) \otimes V(s)^{*}_{1}[-8]$, $K_{1}^{\backsim} = S(V(2)_{0}) \otimes V(s)_{1}$ and zero otherwise, 
with differential
\begin{equation}
\label{eq:k''}
     S(V(2)_{0}) \otimes V(s)^{*}_{1}[-8] \overset{b_{2}^{\backsim}}{\rightarrow} S(V(2)_{0}) \otimes V(s)_{1}     
\end{equation}
of the form
\[     b_{2}^{\backsim} (y \otimes r_{1,a}) = \sum_{i = 1}^{n} \sum_{b = 1}^{s} \Gamma_{a,b}^{i} y x_{i} \otimes z_{b},     \]
where we remind that the action of $x_{i}$ on $y$ is zero for $i = 3, \dots, n$, but we prefer to write the complete sum for convenience. 
The complex $(K_{\bullet}^{\sim},b_{\bullet}^{\sim})$ is similar to the one analyzed in \cite{HS10}, Prop. 3.5, (3.3). 
It is obvious that $H_{0}(K_{\bullet}^{\sim}) \simeq k$, and further $H_{3}(K_{\bullet}^{\sim})=0$, for the map $b_{3}^{\sim}$ is injective. 
Indeed, $b_{3}^{\sim}(y) = 0$, for $y \in S(V(2)_{2})$, means that $y x_{i} = 0$, for $i = 1, 2$, which implies that $y$ vanishes, for $S(V(2)_{0})$ is a domain. 

Moreover, we shall also show that $H_{2}(K_{\bullet}^{\sim})=0$. 
Consider $y = \sum_{i=1}^{n} y_{i} \otimes r_{0,i} \in K_{2}^{\sim}$ in the kernel of $b_{2}^{\sim}$, which is equivalent to the vanishing of 
\[     b_{2}^{\sim}(\sum_{i=1}^{n} y_{i} \otimes r_{0,i}) = \sum_{i,j=1}^{n} y_{i} x_{j}^{2} \otimes x_{i} - y_{i} x_{j} x_{i} \otimes x_{j} 
                                                          = \sum_{i,j=1}^{n} (y_{i} x_{j}^{2} - y_{j} x_{i} x_{j}) \otimes x_{i}.      \]
This means that $\sum_{j=1}^{n} (y_{i} x_{j}^{2} - y_{j} x_{i} x_{j}) = 0$, for all $i = 1, \dots, n$. 
Since the action of $x_{i}$ on $S(V(2)_{0})$ is zero for $i = 3, \dots, n$, we get that $y_{i} (\sum_{j=1}^{2} x_{j}^{2}) = 0$, for all $i = 3, \dots, n$, 
which in turn implies that $y_{i}$ vanishes for $i = 3, \dots, n$. 
Hence, the cycle $y$ has in fact the form $\sum_{i=1}^{2} y_{i} \otimes r_{0,i}$, and satisfies that 
$\sum_{j=1}^{2} (y_{i} x_{j}^{2} - y_{j} x_{i} x_{j}) = 0$, for $i = 1, 2$. 
So, $y$ can be regarded as a cycle of the complex analyzed in \cite{HS10}, Prop. 3.5, (3.3), for $n=2$, whose second homology group 
vanishes. 
By the definition of the differential of this latter complex, we conclude that there exists $x \in S(V(2)_{0})$ such that $y_{i} = x x_{i}$, for $i = 1, 2$, 
which tells us that $b_{3}^{\sim}(x) = y$.  

We shall further prove that $H_{\bullet}(K_{\bullet}^{\backsim})$ also vanishes for $\bullet = 2$, and so for $\bullet \geq 2$. 
Otherwise stated, we will prove that $b_{2}^{\backsim}$ is injective. 
Since $b_{2}^{\backsim}$ is a morphism between two finitely generated free modules (of the same rank) over the commutative domain $S(V_{2})$, 
$b_{2}^{\backsim}$ is injective if and only if its determinant $\det(b_{2}^{\backsim})$ does not vanish. 
Since the matrix representation of $b_{2}^{\backsim}$ in the respective bases $\{r_{1,a}\}_{a = 1, \dots, s}$ and $\{ z_{a} \}_{a = 1, \dots, s}$ is given by 
$(x_{1} \delta_{a,b} + \Gamma_{a,b}^{2} x_{2})_{a,b}$, its determinant has a term of the form $x_{1}^{s}$, and hence it does not vanish. 
As a consequence, $H_{\bullet}(K_{\bullet}^{\backsim})$ vanishes for $\bullet = 2$, and thus $H_{\bullet}(\ym(n,s)^{\Gamma},S(V(2)_{0})) = 0$, for $\bullet \geq 2$. 

It remains to prove that the Hilbert series of $\hat{W}(n,s)^{\Gamma} = H_{1}(\ym(n,s)^{\Gamma},S(V(2)_{0}))$ is the stated above. 
By the previous remarks, we know that $H_{1}(\ym(n,s)^{\Gamma},S(V(2)_{0})) \simeq  H_{1}(K_{\bullet}^{\sim}) \oplus H_{1}(K_{\bullet}^{\backsim})$. 
This follows from the fact that the Euler-Poincar\'e characteristic of the complex 
$(S(V(2)_{0}) \otimes_{\YM(n,s)^{\Gamma}} K_{\bullet}(\ym(n,s)^{\Gamma}), 1 \otimes b_{\bullet}')$ of graded vector spaces, 
with respect to the Euler-Poincar\'e map given by taking Hilbert series, 
coincides with that of its homology, \textit{i.e.} 
\[     1 - \hat{W}(n,s)^{\Gamma}(t) = \frac{1 - n t^{2} - s t^{3} + n t^{6} + s t^{5} - t^{8}}{(1-t^{2})^{2}}.     \]
The Hilbert series of $\hat{W}(n,s)^{\Gamma}$ stated at the beginning follows easily from the previous identity, and the proposition is thus proved. 
\qed

\begin{remark}
If the super Yang-Mills is equivariant, the injectivity of $b_{2}^{\backsim}$ can also be proved as follows. 
Consider the map 
\begin{align*}
   S(V(n)_{0}) \otimes V(s)_{1} &\overset{(b_{2}^{\backsim})^{*}}{\rightarrow} S(V(n)_{0}) \otimes V(s)^{*}_{1}[-8]     
   \\
   z \otimes z_{a} &\mapsto \sum_{i = 1}^{n} \sum_{b = 1}^{s} \tilde{\Gamma}^{i,a,b} z x_{i} \otimes r_{1,b},     
\end{align*}
where we write again all the variables in the previous sum for convenience, even though the action of $x_{i}$ on $S(V(2)_{0})$ vanishes for $i = 3, \dots, n$. 
Condition \eqref{eq:relgammasind} tells us that $((b_{2}^{\backsim})^{*} \circ b_{2}^{\backsim})(z \otimes r_{1,a}) = z.(\sum_{i=1}^{n} x_{i}^{2}) \otimes r_{1,a}$, 
which in turn implies that $b_{2}^{\backsim}$ is injective, for $z.(\sum_{i=1}^{n} x_{i}^{2}) = z.(\sum_{i=1}^{2} x_{i}^{2})$ and 
$S(V(2)_{0})$ is a domain. 
\end{remark}

As a consequence of the previous proposition we obtain the following. 
\begin{theorem}
\label{teo:libre}
Let $n \geq 2$. 
The graded Lie algebra $\hat{\tym}(n,s)^{\Gamma}$ is free with space of generators isomorphic to $\hat{W}(n,s)^{\Gamma}$. 
\end{theorem}
\noindent\textbf{Proof.}
Since $\U(\ym(n,s)^{\Gamma}) \otimes_{\U(\hat{\tym}(n,s)^{\Gamma})} k \simeq S(V(2)_{0})$, Schapiro's Lemma for graded Lie algebras 
(which is proved in the same way as for Lie algebras or group algebras, \textit{cf.} \cite{Wei94}, Lemma 6.3.2) 
implies that $H_{\bullet} (\ym(n,s)^{\Gamma},S(V(2)_{0})) \simeq H_{\bullet}(\hat{\tym}(n,s)^{\Gamma},k)$. 
The previous computation tells us that $H_{\bullet}(\hat{\tym}(n,s)^{\Gamma},k) \simeq \Tor^{\hat{\TYM}(n,s)^{\Gamma}}_{\bullet} (k,k) = 0$, for $\bullet \geq 2$, 
so the minimal projective resolution of the $\hat{\TYM}(n,s)^{\Gamma}$-module $k$, which is of the form 
$(\Tor_{\bullet}^{\hat{\TYM}(n,s)^{\Gamma}}(k,k),d_{\bullet})$, 
has only two nonvanishing components, for $\bullet = 0, 1$. 
Hence $H_{\bullet}(\hat{\tym}(n,s)^{\Gamma},M) = 0$, for $\bullet \geq 2$ and for any $\hat{\tym}(n,s)^{\Gamma}$-module $M$, 
which implies that $\hat{\tym}(n,s)^{\Gamma}$ is a free graded Lie algebra 
(this is proved in the same way as for Lie algebras, \textit{cf.} \cite{Wei94}, Exer. 7.6.3). 
Moreover, Shapiro's Lemma tells us that $\hat{\tym}(n,s)^{\Gamma} = \f(\hat{W}(n,s)^{\Gamma})$, for 
$\hat{W}(n,s)^{\Gamma} = \hat{\tym}(n,s)^{\Gamma}/[\hat{\tym}(n,s)^{\Gamma},\hat{\tym}(n,s)^{\Gamma}] \simeq H_{1}(\hat{\tym}(n,s)^{\Gamma},k) 
\simeq H_{1}(\ym(n,s)^{\Gamma}, S(V(2)_{0}))$ and the theorem follows. 
\qed

\begin{remark}
Since the Lie ideal $\tym(n,s)^{\Gamma}$ is a subalgebra of the free graded Lie algebra $\hat{\tym}(n,s)^{\Gamma}$, it is also a free graded Lie algebra, 
by Shapiro's Lemma. 
As in the proof of the theorem, taking into account that $\U(\ym(n,s)^{\Gamma}) \otimes_{\U(\tym(n,s)^{\Gamma})} k \simeq S(V(n)_{0})$, 
Shapiro's Lemma tells us that $H_{\bullet} (\ym(n,s)^{\Gamma},S(V(n)_{0})) \simeq H_{\bullet}(\tym(n,s)^{\Gamma},k)$. 
Since $\tym(n,s)^{\Gamma}$ is free, the latter homology group should vanish for $\bullet \geq 2$, $H_{0}(\tym(n,s)^{\Gamma},k) \simeq k$, and 
$H_{1}(\tym(n,s)^{\Gamma},k) \simeq H_{1}(\ym(n,s)^{\Gamma},S(V(n)_{0})) \simeq W(n,s)^{\Gamma}$, 
where $W(n,s)^{\Gamma}$ is the graded vector space of generators of $\tym(n,s)^{\Gamma}$. 
Its Hilbert series can be computed using again the fact that the Euler-Poincar\'e characteristic of a complex of graded vector spaces, 
with respect to the Euler-Poincar\'e map given by taking Hilbert series, coincides with that of its homology, applied to the complex 
$(S(V(n)_{0}) \otimes_{\YM(n,s)^{\Gamma}} K_{\bullet}(\ym(n,s)^{\Gamma}), 1 \otimes b_{\bullet}')$. 
We get that 
\[     W(n,s)^{\Gamma}(t) = \frac{(1-t^{2})^{n}-1+n t^{2} + s t^{3} - s t^{5} - n t^{6} +t^{8}}{(1-t^{2})^{n}}.     \]
Note that the graded vector space given by the even part of $W(n,s)^{\Gamma}$ is isomorphic to the graded vector space $W(n)$ considered in \cite{HS10h}, Section 3, 
provided with the special grading. 
\end{remark}

For completeness, we shall also analyze the case of super Yang-Mills algebras $\ym(n,s)^{\Gamma}$ when $n=1$. 
In the case $s = 0$, the super Yang-Mills algebra $\ym(1,0)^{0}$ is just the one-dimensional super Lie algebra concentrated in degree zero. 
We shall now restrict ourselves to the case $s \neq 0$. 
As noted before, the nondegeneracy of $\Gamma$ implies that $\ym(1,s)^{\Gamma} \simeq k.x_{1} \times \h(1,s)$, 
where $\h(1,s) \simeq \f(z_{1},\dots,z_{s})/\cl{\sum_{a=1}^{s} [z_{a},z_{a}]}$. 

The cases $\ym(1,1)$ and $\ym(1,2)$ are nilpotent and finite dimensional, and its representation theory can be analyzed using the below recalled 
Kirillov orbit method.  
In particular, $\ym(1,1)^{\Gamma}$ is a supercommutative super Lie algebra of super dimension $(1,1)$. 
The super Lie algebra $\ym(1,2)^{\Gamma}$ has super dimension $(2,2)$, with basis $x_{1}$, $z_{1}$, $z_{2}$ and $[z_{2},z_{2}]$, where $[z_{1},z_{1}]=-[z_{2},z_{2}]$, and all other brackets vanish. 
As a consequence, since the enveloping algebra of a finite dimensional super Lie algebra is noetherian, we get that $\YM(1,1)^{\Gamma}$ 
and $\YM(1,2)^{\Gamma}$ are noetherian. 

We shall now suppose that $s \geq 3$. 
In this case, it is easy to prove that the Lie ideal $\k(1,n)$ of $\h(1,s)$ given by $\cl{z_{3}, \dots, z_{s}, [z_{1},z_{2}]} + F^{7}\h(1,n)$, where 
$F^{7}\h(1,n)$ is the super vector space formed of elements of $\h(1,s)$ of degree greater than or equal to $9$. 
It is also an ideal when regarded inside $\ym(1,s)^{\Gamma}$. 
Note that $\U(\h(1,n)/\k(1,n)) \simeq k\cl{z_{1}, z_{2}}/\cl{z_{1}^{2}+z_{2}^{2},[z_{1},z_{2}]}$, where the $[z_{1},z_{2}]$ is the supercommutator of 
$z_{1}$ and $z_{2}$. 
Using the Diamond Lemma, it is easy to see that a (homogeneous) basis of it is given by $z_{1}^{\alpha} z_{2}^{\beta}$, 
where $\alpha \in \{ 0 , 1\}$, and $\beta \in \NN_{0}$, and the multiplication is determined by 
\[     z_{2}^{\beta} z_{1} z_{2}^{\beta'} = (-1)^{\beta} z_{1} z_{2}^{\beta + \beta'},     \]
and
\[     z_{1} z_{2}^{\beta} z_{1} z_{2}^{\beta'} = (-1)^{\beta} z_{2}^{\beta + \beta' + 2}.     \]

The following proposition is a direct consequence of the expression of the complex \eqref{eq:hns} and the previous description of $\U(\h(1,n)/\k(1,n))$. 
\begin{proposition}
Let $s \geq 3$. 
The homology groups $H_{\bullet} (\h(1,s), \U(\h(1,n)/\k(1,s)))$ of the graded Lie algebra $\h(1,s)$ with coefficients in 
the module $\U(\h(1,n)/\k(1,n))$, which is obviously a graded right $H(1,s)$-module, are given by 
\[     H_{\bullet} (\h(1,s), \U(\h(1,s)/\k(1,s))) \simeq \begin{cases}
                                                   k, &\text{if $\bullet = 0$,}
                                                   \\
                                                   \tilde{W}(1,s), &\text{if $\bullet = 1$,}
                                                   \\

                                                   0,             &\text{else,}
                                                   \end{cases}
\]
where $\tilde{W}(1,s)$ is the graded vector space with Hilbert series 
\[     \tilde{W}(1,s)(t) = (s-2) t^{3} + (2s-3) t^{6} + \sum_{k \geq 3} (2s - 4) t^{3 k}.     \]

\end{proposition}

The arguments used in Theorem \ref{teo:libre} also yield
\begin{theorem}
\label{teo:libre2}
Let $s \geq 3$. 
The subalgebra $\k(1,s)$ of the graded Lie algebra $\h(1,n)$, and also of $\ym(1,s)^{\Gamma}$, is free with space of generators isomorphic to 
$\tilde{W}(1,s)$. 
\end{theorem}

As a corollary of the previous theorems we obtain the following result. 
\begin{corollary}
\label{cor:nonoeth}
Let $n \geq 2$ or $n = 1$ and $s \geq 3$. 
In either case the super algebra $\YM(n,s)^{\Gamma}$ is not (left nor right) noetherian. 
\end{corollary}
\noindent\textbf{Proof.}
Let us denote by $\h$ a free super Lie algebra inside $\ym(n,s)^{\Gamma}$ such that the dimension of (the underlying vector space of) $\h$ is greater than $1$, 
which exists by the previous theorems. 
The proposition is now proved by an analogous argument as the one presented in the three paragraphs of \cite{HS10} before Remark 3.14. 
Since the extension $\U(\ym(n,s)^{\Gamma}) \supseteq \U(\h)$ is free 
(\textit{i.e.} $\U(\ym(n,s)^{\Gamma})$ is a free (left) module over $\U(\h)$), any (left) ideal $I \subseteq \U(\h)$ satisfies that 
$\U(\ym(n,s)^{\Gamma}).I \cap \U(\h) = I$. 
The fact that the super algebra $\U(\h)$ is not noetherian (because it is free with more than one generator) implies that the super algebra $\U(\ym(n,s)^{\Gamma})$ is not noetherian. 
\qed

Even though the previous algebras are not noetherian, we will prove that they are (left and right) coherent. 
We recall that a (say left) finitely generated module $M$ over graded algebra $A$ is called \emph{coherent} if all its finitely generated submodules are finitely presented. 
It is easy to see that the category of coherent modules over a graded ring is an abelian subcategory of the category of all modules ${}_{A}\Mod$ 
over the graded algebra $A$. 
Moreover, the graded algebra $A$ is called left (resp. right) \emph{coherent} if the category of finitely presented left (resp. right) modules coincides with the category of coherent modules, or otherwise stated if the former is abelian. 
We stress that all modules here are graded. 
Also note that the previous condition can be easily seen to be equivalent to either of the following statements: 
any finitely generated submodule of a finitely presented left (resp. right) module is finitely presented, or any finitely generated left (resp. graded) ideal is finitely presented. 
We shall denote by ${}_{A}\mod$ the category of coherent (graded) modules over the graded algebra $A$. 

Note that the super Yang-Mills algebras are (graded) Hopf algebras, for which the notions of left and right noetherianity coincide, 
and the same applies to coherency. 
So, we shall just refer to these algebras as noetherian, or coherent. 

We first recall the following result due to D. Piontkovski.
\begin{proposition}
\label{prop:piont}
Let $A$ be a graded algebra and $J$ a two-sided ideal of $A$ (different from $A$), which is free as a left module. 
Then, if the quotient graded algebra $A/J$ is right noetherian, $A$ is right coherent. 
\end{proposition}
\noindent\textbf{Proof.}
See \cite{Pi08}, Prop. 3.2. 
\qed

We have the following immediate consequence. 
\begin{corollary}
\label{coro:coh}
Let $\g$ be a $\NN$-graded Lie algebra and $\h$ a Lie ideal of $\g$ ($\h \neq \g$). 
Then, the $\U(\g)$-module $\U(\g)\h$ is free if and only if $\h$ is a free graded Lie algebra. 
As a consequence, assuming that the quotient graded algebra $\U(\g/\h)$ is noetherian and $\h$ is a free graded Lie algebra, $\U(\g)$ is coherent. 
\end{corollary}
\noindent\textbf{Proof.}
Note first that $\U(\g)\h = \h \U(\g)$ is a two-sided ideal, and $\U(\g/\h) \simeq \U(\g)/(\U(\g)\h)$, so the second statement follows directly from the first one and Proposition \ref{prop:piont}. 
In consequence, we only have to prove that $\U(\g)\h$ is a free $\U(\g)$-module if and only if $\h$ is a free graded Lie algebra. 
The former is equivalent to show that $\Tor^{\U(\g)}_{\bullet}(k,\U(\g)\h)$ vanishes for $\bullet \geq 1$, since $\U(\g) \h$ is a bounded below graded $\U(\g)$-module, 
and the latter is equivalent to the vanishing of $\Tor^{\U(\h)}_{\bullet}(k,k) \simeq \Tor^{\U(\g)}_{\bullet}(k,\U(\g/\h))$, for $\bullet \geq 2$, since $\U(\h)$ is a graded connected algebra, 
as explained in Theorem \ref{teo:libre}. 
The previous isomorphism follows from Schapiro's Lemma, since 
$\Tor^{\U(\h)}_{\bullet}(k,k) \simeq H_{\bullet}(\h,k) \simeq H_{\bullet}(\g,\U(\g/\h)) \simeq \Tor^{\U(\g)}_{\bullet}(k,\U(\g/\h))$.  
It thus suffices to prove that both homology groups $\Tor^{\U(\g)}_{\bullet}(k,\U(\g)\h)$ and $\Tor^{\U(\g)}_{\bullet+1}(k,\U(\g/\h))$ are isomorphic, for $\bullet \geq 1$. 

Consider the short exact sequence of $\U(\g)$-modules 
\[     0 \rightarrow \U(\g) \h \rightarrow \U(\g) \rightarrow \U(\g/\h) \rightarrow 0,     \]
which induces a long exact sequence on torsion groups 
\[     \dots \rightarrow \Tor_{i+1}^{\U(\g)}(k,\U(\g/\h)) \rightarrow \Tor_{i}^{\U(\g)}(k,\U(\g)\h) \rightarrow \Tor_{i}^{\U(\g)}(k,\U(\g)) \rightarrow \Tor_{i}^{\U(\g)}(k,\U(\g/\h)) \rightarrow \Tor_{i-1}^{\U(\g)}(k,\U(\g)\h) \rightarrow \dots     \]
The vanishing of $\Tor_{\bullet}^{\U(\g)}(k,\U(\g))$ for $\bullet \geq 1$ (because $\U(\g)$ is free) implies that 
$\Tor_{\bullet+1}^{\U(\g)}(k,\U(\g/\h)) \simeq \Tor_{\bullet}^{\U(\g)}(k,\U(\g)\h))$ for $\bullet \geq 1$, 
which proves the corollary. 
\qed

Theorems \ref{teo:libre} and \ref{teo:libre2} together with the previous corollary now yield the result:
\begin{corollary}
Let $n \geq 2$, or $n = 1$ and $s \geq 3$. 
Then the super Yang-Mills algebra $\YM(n,s)^{\Gamma}$ is (left and right) coherent. 
\end{corollary}

\begin{remark}
\label{rem:coh}
Note that $\ym(0,s)^{0}$ is a free graded algebra on the generators $z_{1}, \dots, z_{s}$, and hence it is obviously seen to be coherent (applying the Corollary \ref{coro:coh} for $\h = \f(z_{2}, \dots, z_{s})$). 
For the other values of the parameters ($n = 1$ and $s = 1$ or $s = 2$), the super Yang-Mills algebra is noetherian, so \textit{a fortiori} coherent. 
As a consequence, we see that the super Yang-Mills algebras are coherent for all values of the parameters $n$ and $s$. 
\end{remark}

Since $\YM(n,s)^{\Gamma}$ is AS-regular with Gorenstein parameter $8$, 
we may consider the (finite dimensional) \emph{Beilinson algebra} $\nabla \YM(n,s)^{\Gamma}$ associated to it, 
\textit{i.e.} 
\[     \nabla \YM(n,s)^{\Gamma} = \begin{pmatrix} 
                                  \YM(n,s)^{\Gamma}_{0} & \YM(n,s)^{\Gamma}_{1} & \dots & \YM(n,s)^{\Gamma}_{7} 
                                  \\ 
                                  0 & \YM(n,s)^{\Gamma}_{0} & \dots & \YM(n,s)^{\Gamma}_{6} 
                                  \\
                                  \vdots & \vdots & \dots & \vdots
                                  \\
                                  0 & 0 & \dots & \YM(n,s)^{\Gamma}_{0}
                                  \end{pmatrix},     \]  
with the obvious matrix multiplication. 
Equivalently, it can also be easily defined by a quiver algebra with relations. 

We now recall for a connected $\NN_{0}$-graded algebra $A$ the definition of the (abelian) quotient categories 
$\mathrm{Tails}(A) = {}_{A}\Mod/{}_{A}\mathrm{Tors}$ and $\mathrm{tails}(A) = {}_{A}\mod/{}_{A}\mathrm{tors}$, 
where ${}_{A}\mathrm{Tors}$ is the category of torsion (graded) modules over $A$, 
\textit{i.e.} the modules $M$ that satisfy that for any $m \in M$, there is $i \in \NN$ such that $A_{\geq i}.m = 0$,  
and ${}_{A}\mathrm{tors}$ is the category of torsion coherent modules over $A$ (\textit{cf.} \cite{AZ94}, Section 1). 

Then, using \cite{MM11}, Thm. 4.12 and 4.14, we obtain the 
\begin{proposition}
Let $n \geq 2$, or $n = 1$ and $s \geq 3$. 
Then, there exists equivalences of triangulated categories 
\begin{align*}
   D(\mathrm{Tails}(\YM(n,s)^{\Gamma})) &\simeq D({}_{\nabla \YM(n,s)^{\Gamma}}\Mod),
   \\
   D^{b}(\mathrm{tails}(\YM(n,s)^{\Gamma})) &\simeq D^{b}({}_{\nabla \YM(n,s)^{\Gamma}}\mod).
\end{align*}
\end{proposition}

As a further consequence of the freeness of the subalgebra $\hat{\tym}(n,s)^{\Gamma}$, we obtain the following result.
\begin{proposition}
Let $n \geq 2$. 
Then, the super Yang-Mills algebra $\YM(n,s)^{\Gamma}$ is a semiprimitive domain, \textit{i.e.}, it has vanishing Jacobson radical and does not have zero divisors.
\end{proposition}
\noindent\textbf{Proof.}
We first prove that the super Yang-Mills algebra is a domain. 
In order to do so, note that, since all odd elements of the super Yang-Mills algebra $\ym(n,s)^{\Gamma}$ are in fact included in the free super Lie algebra $\hat{\tym}(n,s)^{\Gamma}$, 
we have that $[y,y] \neq 0$, for all its odd elements $y$. 
Then, \cite{AL85}, Thm. 2.7, implies that $\YM(n,s)^{\Gamma}$ is a domain. 

We shall now prove that the super Yang-Mills algebra is semiprimitive, \textit{i.e.} the Jacobson radical of $\YM(n,s)^{\Gamma}$ vanishes. 
We recall that the Jacobson radical of a super algebra and of its underlying algebra coincide (see \cite{CM84}, Thm. 4.4, (3)). 
Also note that the Jacobson radical of a free algebra is zero, since for any nonzero element $w$ of a free algebra, $1 - u w$ cannot be invertible for all elements $u$ in the free algebra 
(take $u$ to be nonzero and noninvertible). 

Let $n \geq 2$. 
Then the collection of one-codimensional inclusions of Lie ideals 
$\hat{\tym}(n,s)^{\Gamma} \subseteq \h(n,s) \subseteq \ym(n,s)^{\Gamma}$ tells US that $\U(\h(n,s)) \simeq \U(\hat{\tym}(n,s)^{\Gamma})[x_{2},\delta_{2}]$, 
where $\delta_{2}$ is the derivation on $\U(\hat{\tym}(n,s)^{\Gamma})$ induced by $\ad(x_{2})$, and 
$\U(\ym(n,s)^{\Gamma}) \simeq \U(\h(n,s))[x_{1},\delta_{1}]$, where $\delta_{1}$ is the derivation on $\U(\h(n,s))$ induced by $\ad(x_{1})$. 
As a consequence, the super Yang-Mills algebra $\U(\ym(n,s)^{\Gamma})$ is a sequence of Ore extensions with derivations of a free algebra, for $n \geq 2$. 

Also, note that, if $A \subseteq B$, is an Ore extension, then it is free, and in particular, by the proof of Corollary \ref{cor:nonoeth}, for 
any left ideal $I$ of $A$, we have that $B.I \cap A = I$. 
This in turn implies that $J(B) \cap A \subseteq J(A)$. 
Indeed, consider the collection $\mathcal{S}$ of maximal left ideals $J$ of $B$ such that they contain $B.I$, for some $I$ a maximal left ideal of $A$. 
Note that $B.I \neq B$, since $B.I \cap A = I$. 
Moreover, if $J \in \mathcal{S}$, and $J \supseteq B.I$, for $I$ a maximal ideal of $A$, then $I = B.I \cap A \subseteq  J \cap A \neq A$, and hence 
$J \cap A = I$, for $I$ is maximal. 
We remark that $J \cap A \neq A$, because $1_{A} = 1_{B} \notin J$. 
Let $\mathcal{S}'$ denote the collection of all maximal left ideals of $A$. 
Since $J(B)$ is the intersection of all maximal left ideals of $B$, then 
\[     J(B) \cap A \subseteq \bigcap_{J \in \mathcal{S}} J \cap A = \bigcap_{I \in \mathcal{S}'} I = J(A).     \]

The proposition now follows from \cite{FKM83}, Thm. 3.2. 
\qed

\section{The main result on representations of super Yang-Mills algebra}
\label{sec:main}

The aim of this last section is to prove that (most of) the Clifford-Weyl super algebras $\Cliff_{q}(k) \otimes A_{p}(k)$  
are epimorphic images of all super Yang-Mills algebras $\YM(n,s)^{\Gamma}$, under certain assumptions.
This will rely on our previous study of the Lie ideal $\hat{\tym}(n,s)^{\Gamma}$.
As a consequence, the representations of such Clifford-Weyl super algebras $\Cliff_{q}(k) \otimes A_{p}(k)$ 
are also representations of $\YM(n,s)^{\Gamma}$, which is the analogous result to the one proved in \cite{HS10}. 

On the one hand, since the super algebra $\Cliff_{2q}(k) \simeq M_{2^{q}}(k)$ is Morita equivalent to the super algebra $k$, 
we easily conclude that the Clifford-Weyl super algebra $\Cliff_{2q}(k) \otimes A_{p}(k)$ is Morita equivalent to the super algebra $A_{p}(k)$, and 
$\Cliff_{2q + 1}(k) \otimes A_{p}(k)$ is Morita equivalent to $\Cliff_{1}(k) \otimes A_{p}(k)$. 
Furthermore, the category of representations of the super algebra $\Cliff_{1}(k) \otimes A_{p}(k)$ is equivalent to the category of representations of the Weyl algebra $A_{p}(k)$. 
In both cases we see that the representations we obtain can be understood as induced by those of the Weyl algebra (see Remark \ref{rem:morita}). 

\subsection{Some prerequisites}

We shall now briefly recall a version of the Kirillov orbit method and the Dixmier map for nilpotent super Lie algebras, which we will employ. 
We shall use the conventions of \cite{Her10}, to which we refer for the details and the bibliography therein. 
We recall that a bilateral ideal $I \triangleleft A$ of a super algebra $A$ is called \emph{primitive} if it is the annihilating ideal of a simple $A$-module, 
and it is called \emph{maximal} if $I \neq A$ and if it is maximal in the lattice of bilateral ideals of $A$, ordered by inclusion.
Every maximal ideal is clearly primitive.

We have the following proposition. 
\begin{proposition}
\label{prop:lieweyl}
Let $I \triangleleft \U(\g)$ be a bilateral ideal of the universal enveloping algebra of a nilpotent super Lie algebra $\g$ of finite dimension.
The following are equivalent:
\begin{itemize}
\item[(i)] $I$ is primitive.
\item[(ii)] $I$ is maximal.
\item[(iii)] There exist $p, q \in \NN_{0}$ such that $\U(\g)/I \simeq \Cliff_{p}(k) \otimes A_{q}(k)$, where $\Cliff_{p}(k)$ is the Clifford (super) algebra over $k$ 
and $A_{q}(k)$ is the $q$-th Weyl algebra, which is concentrated in degree zero.
\item[(iv)] $I$ is the kernel of a simple representation of $\U(\g)$.
\end{itemize}
\end{proposition}
\noindent\textbf{Proof.} 
See \cite{Let89}, Prop. 3.3, and \cite{Her10}, Prop. 4.13.
\qed

We remark that the Clifford-Weyl super algebra $\Cliff_{q}(k) \otimes A_{p}(k)$ in the previous theorem has the $\ZZ/2\ZZ$-grading induced by usual grading of the Clifford (super) algebra 
$\Cliff_{q}(k)$ and by considering the Weyl algebra $A_{p}(k)$ to be concentrated in degree zero. 

If $I \triangleleft \U(\g)$ is a bilateral ideal satisfying either of the previous equivalent conditions, the pair of nonnegative integers $(p,q)$ 
(uniquely determined) such that $\U(\g)/I \simeq \Cliff_{p}(k) \otimes A_{q}(k)$ is called the \emph{weight} of the ideal $I$. 

Let us suppose that $\g$ is a nilpotent super Lie algebra of finite dimension. 
Given $f \in \g^{*}_{0}$, a \emph{polarization} of $\g$ at $f$ is a subalgebra $\h \subset \g$ such that it is \emph{subordinated} to $f$, 
\textit{i.e.} $f([\h,\h]) = 0$, 
and it is in fact a maximal subspace of the super vector space underlying $\g$ with respect to the property that 
the bilinear form $f([\place,\place])$ vanishes on it. 
It is easily verified that the super dimensions of all polarizations of $\g$ at $f$ coincide and, furthermore, 
that $\h$ is a polarization of $\g$ at $f$ if and only if $\h$ is a subalgebra subordinated to $f$ whose super dimension coincides 
with the one of a polarization of $\g$ at $f$ (see \cite{Her10}, Subsection 3.4). 

We shall now explain the connection between primitive ideals and even linear functionals for a nilpotent super Lie algebra. 
If $f \in \g^{*}_{0}$ is a linear functional and $\h_{f}$ a polarization at $f$,
we may define a representation of $\h_{f}$ on the vector space $k.v_{f}$ of dimension $1$ by means of 
$x.v_{f} = f(x) v_{f}$, for $x \in \h_{f}$.
Therefore, we can consider the induced $\U(\g)$-module $V_{f} = \U(\g) \otimes_{\U(\h_{f})} k.v_{f}$.
If we denote the corresponding action by $\rho : \U(\g) \rightarrow \End_{k}(V_{f})$, $I(f) = \Ker(\rho)$ is a bilateral ideal of the 
enveloping algebra $\U(\g)$.
In the previous notation we have omitted the polarization in $I(f)$, due to the following proposition, which states even more. 
\begin{proposition}
Let $\g$ be a nilpotent super Lie algebra of finite dimension, let $f \in \g^{*}_{0}$ and let $\h_{f}$ and $\h'_{f}$ be two polarizations of $f$.
If we denote $\rho : \U(\g) \rightarrow \End_{k}(V_{f})$ and $\rho' : \U(\g) \rightarrow \End_{k}(V'_{f})$ 
the corresponding representations constructed following the previous method, then $\Ker(\rho) = \Ker(\rho')$. 
This ideal is primitive. 
 
On the other hand, if $I$ is a primitive ideal of $\U(\g)$, then there exists $f \in \g^{*}_{0}$ such that $I = I(f)$.
\end{proposition}
\noindent\textbf{Proof.} 
See \cite{Her10}, Thm. 4.5, Thm. 4.7 and Thm. 4.9. 
\qed

The weight of a primitive ideal $I(f)$ is given by $(p,q) = (\dim(\g_{0}/\g^{f}_{0})/2, \dim(\g_{1}/\g^{f}_{1}))$, 
where $\g^{f} = (\g^{f}_{0}, \g^{f}_{1})$ is the kernel of the superantisymmetric bilinear form $f([\place,\place])$ determined by $f$ on $\g$ 
(see \cite{Her10}, Prop. 4.13). 

The group $\Aut(\g_{0})$ is an algebraic group whose associated Lie algebra is $\Der(\g_{0})$.
Since the super Lie algebra $\g$ is nilpotent, the Lie algebra given by its even part $\g_{0}$ is also so, and Lie algebra given by the ideal of inner derivations $\InnDer(\g_{0})$ in $\Der(\g_{0})$ is algebraic. 
The irreducible algebraic group $\mathcal{A}d_{0}$ associated to $\InnDer(\g_{0})$ is called the \emph{adjoint (algebraic) group} of $\g_{0}$.
It is a subgroup of $\Aut(\g_{0})$.
As a consequence, the group $\mathcal{A}d_{0}$ acts on the Lie algebra $\g_{0}$, so it also acts on $\g^{*}_{0}$ with the dual action, which is called \emph{coadjoint}. 
\begin{proposition}
Let $\g$ be a nilpotent super Lie algebra of finite dimension and let $f$ and $f'$ be two even linear forms on $\g$, \textit{i.e.} $f, f' \in \g^{*}_{0}$. 
If $I(f)$ and $I(f')$ are the corresponding bilateral ideals of $\U(\g)$, then $I(f) = I(f')$ if and only if there is $g \in \mathcal{A}d_{0}$ such that 
$f = g.f'$.
\end{proposition}
\noindent\textbf{Proof.} 
See \cite{Her10}, Prop. 4.12.
\qed

The previous results imply that, for a nilpotent super Lie algebra of finite dimension there exists an explicit bijection
\[     I : \g^{*}_{0}/\mathcal{A}d_{0} \rightarrow \Prim (\U(\g))     \]
between the set of equivalence classes of even linear forms on $\g$ under the coadjoint action and the set of primitive ideals of the super algebra $\U(\g)$.

We also recall that, if $\g$ be a finite dimensional nilpotent super Lie algebra and $\h$ a Lie ideal, 
given $I \triangleleft \U(\h) \subset \U(\g)$ a two sided ideal in the enveloping algebra of $\h$, one defines the \emph{stabilizer} of the ideal $I$ in $\g$ to be  
\[     \mathfrak{st}(I, \g) = \{ x \in \g : [x , I] \subset I \}.     \]
Let us further suppose that there exists $f \in \g_{0}^{*}$ such that $I = I(f|_{\h_{0}})$. 
By \cite{Her10}, Prop. 4.16, if $\g' = \{ x \in \g : f([x,\h])=0 \}$, then 
\[     \mathfrak{st}(I, \g) \supseteq \g' + \h.     \]

\subsection{Main theorem}

We have first the following result: 
\begin{proposition}
\label{teo:yangmillsweyl}
Let $n, s, p \in \NN$ be positive integers, satisfying that $n \geq 3$. 
There exists a surjective homomorphism of super algebras
\[     \U(\ym(n,s)^{\Gamma}) \twoheadrightarrow A_{p}(k).     \]
Furthermore, there exists $l \in \NN$ such that we can choose this homomorphism satisfying that it factors through the quotient
$\U(\ym(n,s)^{\Gamma}/F^{l}\ym(n,s)^{\Gamma})$
\[
\xymatrix@C-20pt@R-10pt { \U(\ym(n,s)^{\Gamma}) \ar@{->>}[rr] \ar@{->>}[rd] &
& A_{p}(k)
\\
& \U(\ym(n,s)^{\Gamma}/F^{l}\ym(n,s)^{\Gamma}) \ar@{->>}[ru] & }
\]
\end{proposition}
\noindent\textbf{Proof.}
This is a direct consequence of the surjective morphism of graded algebras given by $\YM(n,s)^{\Gamma} \rightarrow \YM(n)$ described at the end of 
Subsection \ref{subsec:defi} and \cite{HS10}, Coro. 4.5. 
\qed

On the other hand, we also have the:
\begin{theorem}
\label{teo:yangmillsclifweyl}
Let $n, s, p, q \in \NN$ be positive integers, satisfying that $n \geq 3$ and $s \geq 1$. 
We suppose further that either $p \geq 3$, or $p = 2$ and $q \geq 2$. 
Then, there exists a surjective homomorphism of super algebras
\[       \U(\ym(n,s)^{\Gamma}) \twoheadrightarrow \Cliff_{q}(k) \otimes A_{p}(k).     \]
Furthermore, there exists $l \in \NN$ such that we can choose this morphism in such a way that it factors through the quotient
$\U(\ym(n,s)^{\Gamma}/F^{l}(\ym(n,s)^{\Gamma}))$
\[
\xymatrix@C-20pt@R-10pt { \U(\ym(n,s)^{\Gamma}) \ar@{->>}[rr] \ar@{->>}[rd] &
& \Cliff_{q}(k) \otimes A_{p}(k)
\\
& \U(\ym(n,s)^{\Gamma}/F^{l}\ym(n,s)^{\Gamma}) \ar@{->>}[ru] & }
\]
\end{theorem}
\noindent\textbf{Proof.}
We know that $\ym(n,s)^{\Gamma} = V(2)_{0} \oplus \hat{\tym}(n,s)^{\Gamma}$ as graded vector spaces.
As we have proved in Theorem \ref{teo:libre}, the Lie ideal $\hat{\tym}(n,s)^{\Gamma}$, considered as a graded Lie algebra, is a free graded Lie algebra 
generated by a graded vector space $\hat{W}(n,s)^{\Gamma}$, that is,
\[     \hat{\tym}(n,s)^{\Gamma} \simeq \f(\hat{W}(n,s)^{\Gamma}).     \]
We also point out that, by the computation of the Hilbert series of $\hat{W}(n,s)^{\Gamma}$, this space has an infinite number of 
linearly independent even and odd elements (for $n \geq 3$ and $s \in \NN$). 

We introduce the following notation, that we will need in the proof. 
Given a super vector space $X$ with a nondegenerate (even) superantisymmetric bilinear form $\Omega$, consider the super vector space 
$\mathfrak{heis}(X) = X \oplus k.z$. 
It is a super Lie algebra provided with the bracket given by declaring that $z$ is central and $[x,x'] = \Omega(x,x') z$, 
and it is called the \emph{Heisenberg super Lie algebra} defined by $(X,\Omega)$. 
It is uniquely determined by the super dimension $(d,d')$ of $X$, so it will also be denoted by $\mathfrak{heis}_{d,d'}$. 
More concretely, if the super dimension of $X$ is $(2r, 2t')$, $\mathfrak{heis}_{2r, 2t'}$ is the super Lie algebra with basis given by even elements 
$q_{1}, \dots, q_{r}, p_{1}, \dots, p_{r}, z$ and odd ones $a_{1}, \dots, a_{t'}, b_{1}, \dots, b_{t'}$, such that $z$ is central, 
$[q_{i},p_{j}] = \delta_{i,j} z$, $[a_{i},b_{j}] = \delta_{i,j} z$ and the other brackets vanish. 
If the super dimension of $X$ is $(2r, 2t'+1)$, $\mathfrak{heis}_{2r, 2t'+1}$ is the super Lie algebra with basis given by even elements 
$q_{1}, \dots, q_{r}, p_{1}, \dots, p_{r}, z$ and odd ones $a_{1}, \dots, a_{t'}, b_{1}, \dots, b_{t'},c$, such that $z$ is central, 
$[q_{i},p_{j}] = \delta_{i,j} z$, $[a_{i},b_{j}] = \delta_{i,j} z$, $[c,c] = z$ and all other brackets vanish. 
It is easy to see that $\cl{z-1} \subseteq \U(\mathfrak{heis}_{2r,t})$ is primitive, and in fact 
$\U(\mathfrak{heis}_{2r,t})/\cl{z-1} \simeq \Cliff_{t}(k) \otimes A_{r}(k)$ (see \cite{BM89}, 0.2, (a) and (b)). 
Denote the projection induced by the previous quotient by $\Pi_{\mathfrak{heis}_{2r,t}}$. 
Hence, it suffices to prove the theorem for $\U(\mathfrak{heis}_{2p,q})$ instead of $\Cliff_{q}(k) \otimes A_{p}(k)$. 

A morphism of super Lie algebras from $\hat{\tym}(n,s)^{\Gamma}$ to $\mathfrak{heis}_{2r,t}$ induces a respective morphism of super algebras 
from $\U(\hat{\tym}(n,s)^{\Gamma})$ to $\U(\mathfrak{heis}_{2r,t})$. 
Since $\hat{\tym}(n,s)^{\Gamma}$ is free with space of generators $\hat{W}(n,s)^{\Gamma}$, the former is equivalent to give a morphism of 
super vector spaces from $\hat{W}(n,s)^{\Gamma}$ to $\mathfrak{heis}_{2r,t}$.
The morphism of super algebras will be surjective if the image of the corresponding morphism of super vector spaces is also so. 
Since $\mathfrak{heis}_{2r,t}$ has finite super dimension, but $\hat{W}(n,s)^{\Gamma}$ has an infinite number of even and odd generators, this can be easily done. 

From now on, we shall exclusively work in the super case (for algebras).
However, we will also keep track of the $\NN$-grading of the super Yang-Mills algebra $\ym(n,s)^{\Gamma}$, its ideal $\hat{\tym}(n,s)^{\Gamma}$ 
and the generator space $\hat{W}(n,s)^{\Gamma}$.

We shall now suppose that either $r \geq 1$, or $r = 0$ and $t \geq 2$, and show that under any of these assumptions, there exists a 
$k$-linear homogeneous morphism (of super vector spaces) of degree $0$ 
\[     \phi : \hat{W}(n,s)^{\Gamma} \rightarrow \mathfrak{heis}_{2r,t}     \]
such that there exists a set of linearly independent homogeneous elements of $\hat{W}(n,s)^{\Gamma}$ 
which are mapped onto a set of homogeneous generators of the Heisenberg super Lie algebra, respecting the degree. 

$x_{3}$ is mapped to zero, $x_{13}$ and $x_{23}$ are mapped to a linearly independent set of two even generators of $\mathfrak{heis}_{2r,t}$, 
and in fact 

If $r \geq 1$, we set $\phi$ such that $\phi(x_{3}) = 0$, $\phi(x_{13}) = p_{1}$ and $\phi(x_{23}) = q_{1}$. 
We also assume that for the set of other basis elements of even parity of the Heisenberg super Lie algebra $z$, $p_{i}$ and $q_{i}$, for $2 \leq i \leq r$, 
there exist a linearly independent set of even homogeneous elements of degree greater than or equal to $6$, 
$w_{i} \in \hat{W}(n,s)^{\Gamma}$ whose image under $\phi$ give the previous elements. 
Let $d_{i}$ be the degree of $w_{i}$, for each ($2 \leq i \leq r$).
Let $j$ be the maximum between $4$ and the degrees $d_{i}$.
We moreover assume that there also exist linearly independent even (resp. odd) homogeneous elements of $\hat{W}(n,s)^{\Gamma}$ of degree greater than $j$ 
mapped onto a set of even (resp. odd) generators of the super Lie algebra $\mathfrak{heis}_{2r,t}$. 
Let us called the maximum of these degrees by $d'$. 
These conditions are easily verified, taking into account that $\hat{W}(n,s)^{\Gamma}$ has infinite dimensional even and odd components.

On the other hand, if $r = 0$ and $t \geq 2$, consider that $\phi(z_{1}) = 0$, $\phi([x_{1},z_{1}]) = a_{1}$ and $\phi([x_{2},z_{1}]) = b_{1}$, 
where $a_{1}, b_{1} \in \mathfrak{heis}_{2r,t}$ are two linearly independent elements. 
We also assume that for the set of other basis elements of the Heisenberg super Lie algebra there exist a linearly independent 
set of even homogeneous elements of degree greater than or equal to $6$, $w_{i} \in \hat{W}(n,s)^{\Gamma}$, 
such that each of the basis elements of the Heisenberg super Lie algebra is the image under $\phi$ of a respective element $w_{i}$ of the same degree. 
Let $d_{i}$ be the degree of $w_{i}$, and let $j$ be the maximum between $4$ and the degrees $d_{i}$.
Again, these conditions are easily verified, taking into account that $\hat{W}(n,s)^{\Gamma}$ has infinite dimensional even and odd components. 

Note that in either case there are a lot of choices for this morphism $\phi$. 
Set $l$ to be $2 d' + 1$.  

The map $\phi$ induces a unique surjective homomorphism $\Phi : \U(\hat{\tym}(n,s)^{\Gamma}) \twoheadrightarrow \U(\mathfrak{heis}_{2r,t})$, 
equivalent to the homomorphism of super Lie algebras
\[     \hat{\tym}(n,s)^{\Gamma} \rightarrow \mathfrak{heis}_{2r,t}.     \]
Since $\mathfrak{heis}_{2r,t}$ is nilpotent, the last morphism may be factorized in the following way
\[     \hat{\tym}(n,s)^{\Gamma} \rightarrow \hat{\tym}(n,s)^{\Gamma}/F^{l}\ym(n,s)^{\Gamma} \rightarrow \mathfrak{heis}_{2r,t},     \]
where the first morphism is the canonical projection.
Hence, the map $\Phi$ may be factorized as 
\[     \U(\hat{\tym}(n,s)^{\Gamma}) \twoheadrightarrow \U(\hat{\tym}(n,s)^{\Gamma}/F^{l}\ym(n,s)^{\Gamma}) \twoheadrightarrow \U(\mathfrak{heis}_{2r,t}).     \]
We have thus obtained a surjective homomorphism of super algebras
\[     \Psi : \U(\hat{\tym}(n,s)^{\Gamma}/F^{l}\ym(n,s)^{\Gamma}) \twoheadrightarrow \U(\mathfrak{heis}_{2r,t}),     \]
where the super Lie algebra $\hat{\tym}(n,s)^{\Gamma}/F^{l}\ym(n,s)^{\Gamma}$ is obviously nilpotent and finite dimensional.
Moreover, it is a Lie ideal of the nilpotent super Lie algebra $\ym(n,s)^{\Gamma}/F^{l}\ym(n,s)^{\Gamma}$.
We have, as graded vector spaces, 
\[     \ym(n,s)^{\Gamma}/F^{l}\ym(n,s)^{\Gamma} = V(2)_{0} \oplus \hat{\tym}(n,s)^{\Gamma}/F^{l}\ym(n,s)^{\Gamma}.     \]

Let $I$ be the kernel of $\Pi_{\mathfrak{heis}_{2r,t}} \circ \Psi$ in $\U(\hat{\tym}(n,s)^{\Gamma}/F^{l}\ym(n,s)^{\Gamma})$.
Taking into account that the quotient of $\U(\hat{\tym}(n,s)^{\Gamma}/F^{l}\ym(n,s)^{\Gamma})$ by $I$ is a Clifford-Weyl super algebra which is simple, $I$ is a maximal two-sided ideal, and then, 
there exists an even linear functional
$f \in (\hat{\tym}(n,s)^{\Gamma}/F^{l}\ym(n,s)^{\Gamma})^{*}$ such that $I = I(f)$.
We fix a polarization $\h_{f}$ of $\hat{\tym}(n,s)^{\Gamma}/F^{l}\ym(n,s)^{\Gamma}$ at $f$.
Let $\bar{f} \in (\ym(n,s)^{\Gamma}/F^{l}\ym(n,s)^{\Gamma})^{*}$ be any (even) extension of $f$. 

Since $\hat{\tym}(n,s)^{\Gamma}/F^{l}\ym(n,s)^{\Gamma}$ is an ideal of the nilpotent super Lie algebra 
$\ym(n,s)^{\Gamma}/F^{l}\ym(n,s)^{\Gamma}$, by \cite{Her10}, Prop. 4.16, we have that the stabilizer 
$\mathfrak{st}(I(f), \ym(n,s)^{\Gamma}/F^{l}\ym(n,s)^{\Gamma})$ includes 
\[     \hat{\tym}(n,s)^{\Gamma}/F^{l}\ym(n,s)^{\Gamma} + (\ym(n,s)^{\Gamma}/F^{l}\ym(n,s)^{\Gamma})',     \] 
and we recall that $(\ym(n,s)^{\Gamma}/F^{l}\ym(n,s)^{\Gamma})'$ is given by 
\[     \{ x \in \ym(n,s)^{\Gamma}/F^{l}\ym(n,s)^{\Gamma} : f([x , \hat{\tym}(n,s)^{\Gamma}/F^{l}\ym(n,s)^{\Gamma}]) = 0 \}.     \]

If $r \geq 1$, we get immediately that $\bar{x}_{3} \in I$, but $\bar{x}_{13}$ and $\bar{x}_{23}$ do not belong to $I$, 
since $\Psi(\bar{x}_{13})$ and $\Psi(\bar{x}_{23})$ do not vanish (they are in fact linearly independent).
Analogously, if $r = 0$ and $t \geq 2$, $\bar{z}_{1} \in I$, but $[\bar{x}_{1},\bar{z}_{1}]$ and $[\bar{x}_{2},\bar{z}_{1}]$ do not belong to $I$, 
since $\Psi([\bar{x}_{1},\bar{z}_{1}])$ and $\Psi([\bar{x}_{2},\bar{z}_{1}])$ are linearly independent.

Let $x \in \ym(n,s)^{\Gamma}/F^{l}\ym(n,s)^{\Gamma}$, then $x = x' + y$, where
\[     x' = \sum_{i=1}^{2} c_{i} \bar{x}_{i} \in V(2)_{0},     \]
and $y \in \hat{\tym}(n,s)^{\Gamma}/F^{l}\ym(n,s)^{\Gamma}$.
Since $[y,I(f)] \subset I(f)$, this implies that $x \in \mathfrak{st}(I(f), \ym(n,s)^{\Gamma}/F^{l}\ym(n,s)^{\Gamma})$ if and only if
$[x' , I(f)] \subset I(f)$.
In particular, if $r \geq 1$,
\[     [x' , \bar{x}_{3}] = \sum_{i=1}^{2} c_{i} [\bar{x}_{i} , \bar{x} _{3}].     \]
If $[x' , \bar{x}_{3}] \in I$, then $\Psi([x' , \bar{x}_{3}]) = 0$, or,
\[     \sum_{i=1}^{2} c_{i} \Psi(\bar{x}_{i3}) = 0,     \]
but since $\Psi(\bar{x}_{13})$ and $\Psi(\bar{x}_{23})$ are linearly independent, we get that $c_{i} = 0$, for all $i = 1, 2$, which gives $x' = 0$. 
In a similar way, if $r = 0$ and $t \geq 2$, 
\[     [x' , \bar{z}_{1}] = \sum_{i=1}^{2} c_{i} [\bar{x}_{i} , \bar{z} _{1}],     \]
so the assumption that $[x' , \bar{z}_{1}] \in I$, implies that 
\[     \sum_{i=1}^{2} c_{i} \Psi([\bar{x}_{i} , \bar{z} _{1}]) = 0,     \]
but since $\Psi([\bar{x}_{1} , \bar{z} _{1}])$ and $\Psi([\bar{x}_{2} , \bar{z} _{1}])$ are linearly independent, we get that $c_{i} = 0$, for all $i = 1, 2$, which again tells us that $x' = 0$.
In any case, we get that $\mathfrak{st}(I(f), \ym(n,s)^{\Gamma}/F^{l}\ym(n,s)^{\Gamma}) = \hat{\tym}(n,s)^{\Gamma}/F^{l}\ym(n,s)^{\Gamma}$.
So, the space $(\ym(n,s)^{\Gamma}/F^{l}\ym(n,s)^{\Gamma})'$ is included in the quotient $\hat{\tym}(n,s)^{\Gamma}/F^{l}\ym(n,s)^{\Gamma}$, 
and, by definition, $\h_{f}$ is also a polarization of $\ym(n,s)^{\Gamma}/F^{l}\ym(n,s)^{\Gamma}$ at $f$. 
In particular, the weight of the ideal $I(\bar{f})$ can be computed easily using \cite{Her10}, Prop. 4.13, previously mentioned, which gives 
$(r+2,t)$. 
Therefore, the quotient of $\YM(n,s)^{\Gamma}$ by the inverse image of $I(\bar{f})$ under the projection 
$\U(\ym(n,s)^{\Gamma}) \rightarrow \U(\ym(n,s)^{\Gamma}/F^{l}\ym(n,s)^{\Gamma})$ is isomorphic to $\Cliff_{t}(k) \otimes A_{r+2}(k)$ and the theorem follows.
\qed

\begin{remark}
We may also study the simplest case cases $\ym(1,1)^{\Gamma}$ and $\ym(1,2)^{\Gamma}$, which are nilpotent and finite dimensional, 
as explained in Paragraph \ref{subsec:tym}. 
Their representation theory can be thus directly analyzed using the Kirillov orbit method.  
In particular, since $\ym(1,1)^{\Gamma}$ is a supercommutative, the only irreducible representations are one-dimensional. 
Concerning the super Lie algebra $\ym(1,2)^{\Gamma}$, it is easy to check that for any even functional $f \in (\ym(1,2)^{\Gamma}_{0})^{*}$, 
we have that $(\ym(1,2)^{\Gamma})^{f} = \ym(1,2)^{\Gamma}$ if $f = 0$, and $(\ym(1,2)^{\Gamma})^{f} = \ym(1,2)^{\Gamma}_{0}$, else. 
This implies that the simple quotients of $\YM(1,2)^{\Gamma}$ are isomorphic to either $k$ or $\Cliff_{2}(k)$. 
\end{remark}

\begin{remark}
\label{rem:otrosn}
Using a similar proof to the one given before, one can also show that, if $n = 1$ and $s \geq 3$, or $n = 2$, then there exists an infinite set of indices 
$(p,q)$ such that  
\[       \U(\ym(n,s)^{\Gamma}) \twoheadrightarrow \Cliff_{q}(k) \otimes A_{p}(k),     \]
where $p \geq 2$ and $q \geq 4$, for $n = 1$ and $s \geq 3$; and $p \geq 4$ and $q \geq 3$, for $n \geq 2$.
Furthermore, there exists $l \in \NN$ such that we can select this morphism in such a way that it factors through the quotient
$\U(\ym(n,s)^{\Gamma}/F^{l}(\ym(n,s)^{\Gamma}))$. 
The previous set can be chosen in such a way that if $(p,q)$ is not in it, then both $(p,q+1)$ and $(p,q-1)$ are. 
\end{remark}

\begin{remark}
\label{rem:morita}
We would like to make some comments on the (abelian) category of representations of the Clifford-Weyl super algebras. 
We remark that the category of representations of a super algebra $A$ is provided with homogeneous $A$-linear morphisms of degree zero. 
Note that, even if a super algebra $A$ is concentrated in degree zero, its category of representations, as a super algebra, does not coincide with the category of representations of the underlying algebra of $A$. 
Indeed, we remark that a representation of a super algebra $A$ concentrated in degree zero is given by a direct sum $M_{0} \oplus M_{1}$ of two modules $M_{0}$ and $M_{1}$ over the underlying algebra of $A$, 
and a morphism from $M_{0} \oplus M_{1}$ to $N_{0} \oplus N_{1}$ is given by a pair $(f_{0},f_{1})$, where $f_{i} : M_{i} \rightarrow N_{i}$ is a morphism of modules over the underlying algebra of $A$, 
for all $i \in \ZZ/2\ZZ$. 
Two super algebras $A$ and $B$ are called \emph{Morita equivalent} if there is an equivalence between the categories of representations of the super algebra $A$ and the super algebra $B$, which commutes with the shift. 
We remark that if two super algebras are Morita equivalent, then the underlying algebras are also Morita equivalent
(\textit{cf.} \cite{GG82}, Section 5, where the authors work with graded algebras, but the statements are analogous, in particular: Lemma 5.1, Cor. 5.2, Prop. 5.3 and Thm. 5.4). 

We now note that the super algebra $\Cliff_{2q}(k) \simeq M_{2^{q}}(k)$ is Morita equivalent to the super algebra $k$. 
This implies that the Clifford-Weyl super algebra $\Cliff_{2q}(k) \otimes A_{p}(k)$ is Morita equivalent to the Weyl algebra $A_{p}(k)$, 
which is regarded as a super algebra concentrated in degree zero, and the super algebras $\Cliff_{2q + 1}(k) \otimes A_{p}(k)$ and 
$\Cliff_{1}(k) \otimes A_{p}(k)$ are also Morita equivalent. 
Moreover, a direct inspection tells us that the category of representations of the super algebra $\Cliff_{1}(k) \otimes A_{p}(k)$ is equivalent to the category of representations of the algebra $A_{p}(k)$. 
We see thus that in either case the representations that we obtain can be understood as induced by those of the Weyl algebras. 
We would also like to point out that several families of representations of the Weyl algebras have been previously studied by Bavula and Bekkert in \cite{BB00}, and, by the theorem, 
they can also be also used to induce representations of the super Yang-Mills algebras. 
\end{remark}



\noindent Estanislao Herscovich,
\\
Fakult\"at f\"ur Mathematik,
\\
Universit\"at Bielefeld,
\\
D-33615 Bielefeld,
\\
Germany,
\\
\href{mailto:eherscov@math.uni-bielefeld.de}{eherscov@math.uni-bielefeld.de}

\end{document}